\documentclass[a4paper]{article}
\usepackage{mathrsfs}
\usepackage{indentfirst,latexsym,bm,amsfonts,amssymb,amsmath,amsfonts,mathrsfs}
\newtheorem{Definition}{Definition}[section]
\newtheorem{Theorem}{Theorem}[section]
\newtheorem{Example}{Example}[section]
\newtheorem{Lemma}{Lemma}[section]
\newtheorem{Corollary}{Corollary}[section]

\title{\bf Demi-linear Analysis III\linebreak
---Demi-distributions with Compact Support}

\author{Li Ronglu\\  
{\it\normalsize Dept. Math., Harbin Institute of Technology, Harbin,
150001, China}\\[3pt]
Zhong Shuhui$^{\ast}$\\
{\it\normalsize Dept. Math., Tianjin University, Tianjin,
300072, China}\\{\it\normalsize  E-mail: shuhuizhong@126.com}\\[3pt]
{\it\normalsize Ctr. Appl. Math., Tianjin University, Tianjin, 300072, China}\\{\it\normalsize  E-mail: shuhuizhong@126.com}\\[3pt]
Kim Dohan\\  
{\it\normalsize Dept. Math., Seoul National University, Seoul,
151-742, Korea}
\\Wu Junde\\ 
{\it\normalsize Dept. Math., Zhejiang University, Hangzhou 310027, China}\\{\it\normalsize  E-mail: wjd@zju.edu.cn}\\[3pt]
}
\begin{document}
\setcounter{footnote}{1}
\renewcommand{\thefootnote}{\fnsymbol{footnote}}
\footnotetext{
{This work was supported by the National Natural Science
Foundation of China (Grant No. 11126165) and the Seed Foundation of
Tianjin University (Grant No. 60302051). }}
\date{}
\maketitle
\def\abstractname{}
\begin{abstract}
\noindent{\bf Abstract. } A series of detailed quantitative results
is established for the family of demi-distributions which is a large
extension of the family of usual distributions.\vspace{6pt}


\end{abstract}\vspace{20pt}

In [1] we show that there is an entirely original generalization of
the basic theory of usual distributions.

As was shown in [1] and [2], the family of demi-distributions is a
large extension of the family of usual distributions, that is, the
family of demi-distributions includes nonlinear functionals as many
as usual distributions, at least.

The theory of demi-distributions not only contains the theory of
usual distributions as a special case but causes a series of
essential changes in the distribution theory. For instance, in the
case of usual distributions the constant distributions are only
solutions of the equation $y'=0$ but in the case of
demi-distributions the equation $y'=0$ has extremely many solutions
which are nonlinear functionals, and every constant is of course a
solution of $y'=0$ [1, Th. 2.3]. Moreover, the family of
demi-distributions is closed with respect to extremely many of
nonlinear transformations such as $|f(\cdot)|$, $|f(\cdot)|^{2/3}$,
$\sin|f(\cdot)|$, $e^{|f(\cdot)|-1}$, etc.

In this paper we carry out a detailed quantitative analysis for
demi-distributions. Our vivid quantitative results show that the
demi-linear mapping introduced in [2] is a very successful
scientific object and, indeed, since the basic principles such as
the equicontinuity theorem and the uniform boundedness principle
hold for the family of demi-linear mappings [2, Th. 3.1, Th. 3.2,
Th. 3.3 , Th. 4.1] and a nice duality theory has established for
demi-linear dual pairs [3, Th. 3.4, Th. 3.12, Th. 3.14, Th. 3.22,
Th. 3.24], it is trivial that the family of demi-linear mappings is
a very valuable extension of the family of linear operators.

\section{Introduction}
Fix an $n\in\mathbb N$ and a nonempty open set $\Omega\subset\mathbb
R^n$. Let
$$C^\infty(\Omega)=\big\{\xi\in\mathbb C^\Omega:\xi\mbox{ is infinitely differentiable in }\Omega\big\},$$
$$C^\infty_0(\Omega)=\big\{\xi\in C^\infty(\Omega):supp\,\xi\mbox{ is compact}\big\},$$
where $supp\,\xi=\overline{\{{x\in\Omega:\xi(x)\neq0}\}}\cap\Omega$
for every $\xi\in\mathbb C^\Omega$ and so if $\xi\in
C^\infty_0(\Omega)$ then the compact
$supp\,\xi=\overline{\{{x\in\Omega:\xi(x)\neq0}\}}\cap\Omega
=\overline{\{x\in\Omega:\xi(x)\neq0\}}\subset\Omega$. For every
$M\subset\Omega$, $C^\infty(M)=\{\xi\in
C^\infty(\Omega):supp\,\xi\subset M\}$ and $C^\infty_0(M)=\{\xi\in
C^\infty_0(\Omega):supp\,\xi\subset M\}$ [4, p.14].

Let $K$ be a compact subset of $\Omega$, that is, $K$ is bounded and
closed in $\mathbb R^n$ and $K\subset\Omega$, and
$k\in\{0,1,2,3,\cdots\}$. Then
$$\|\xi\|_{K,k}=\sum_{|\alpha|\leq k}\sup_K|\partial^\alpha\xi|,\ \xi\in C^\infty(\Omega)$$
defines a seminorm on $C^\infty(\Omega)$, and the family
$\big\{\|\cdot\|_{K,k}:K\mbox{ is compact, }K\subset\Omega,\
k\in\{0\}\cup\mathbb N\big\}$ gives a locally convex Fr\'echet
topology for $C^\infty(\Omega)$, and $C^\infty(\Omega)$ has the
Montel property, i.e., bounded sets in $C^\infty(\Omega)$ are
relatively compact [5, 2.1].

For a compact $K\subset\Omega$ the sequence
$\big\{\|\cdot\|_{K,k}\big\}^\infty_{k=0}$ gives a locally convex
Fr\'echet topology for $C^\infty(K)$. Since
$\Omega=\bigcup^\infty_{j=1}K_j$ where each $K_j$ is compact and
$K_1\subset K_2\subset\cdots$, with the inductive topology using the
inclusion maps, $C^\infty_0(\Omega)=\bigcup_{j=1}^\infty
C^\infty(K_j)$ is a (LF) space which are both barrelled and
bornological. Then $C_0^\infty(\Omega)$ also has the montel
property, and the inclusion map $I:C_0^\infty(\Omega)\rightarrow
C^\infty(\Omega)$ is continuous [5, 2.2].

A distribution $f$ in $\Omega$ is a continuous linear functional on
$C_0^\infty(\Omega)$, that is,
$f:C_0^\infty(\Omega)\rightarrow\mathbb C$ is linear and for every
compact $K\subset\Omega$ there exist $C>0$ and $k\in\{0\}\cup\mathbb
N$ such that
$$\big|f(\xi)\big|\leq C\sum_{|\alpha|\leq k}\sup_K\big|\partial^\alpha\xi\big|,
\ \ \xi\in C_0^\infty(K)=C^\infty(K)\leqno(1.1)$$ [4, Def. 2.1.1,
Th. 2.1.4].

Let $C(0)=\big\{\gamma\in\mathbb C^{\mathbb
C}:\lim_{t\rightarrow0}\gamma(t)=\gamma(0)=0,\
|\gamma(t)|\geq|t|\mbox{ if }|t|\leq1\big\}$. For a topological
vector space $X$, $\mathcal N(X)$ denotes the family of
neighborhoods of $0\in X$.

\begin{Definition}
([2,\ Def. 2.1]) Let $X,\ Y$ be topological vector spaces over the
scalar field $\mathbb K$. A mapping $f:X\rightarrow Y$ is said to be
demi-linear if $f(0)=0$ and there exist $\gamma\in C(0)$ and
$U\in\mathcal N(X)$ such that every $x\in X$, $u\in U$ and
$t\in\big\{t\in\mathbb K:|t|\leq1\big\}$ yield $r,\ s\in\mathbb K$
for which $|r-1|\leq|\gamma(t)|$, $|s|\leq|\gamma(t)|$ and
$f(x+tu)=rf(x)+sf(u)$.
\end{Definition}

Let $\mathscr L_{\gamma,U}(X,Y)$ be the family of demi-linear
mappings related to $\gamma\in C(0)$ and $U\in\mathcal N(X)$, and
let
\begin{align*}
\mathscr K_{\gamma,U}(X,Y)=\big\{f\in&\mathscr
L_{\gamma,U}(X,Y):\mbox{if }x\in X,\ u\in U\mbox{ and
}|t|\leq1,\mbox{ then }\\
&f(x+tu)=f(x)+sf(u)\mbox{ for some }s\mbox{ with
}|s|\leq|\gamma(t)|\big\}.
\end{align*}

If $\gamma(t)=Mt$ with $M\geq1$, then we write that $\mathscr
L_{\gamma,U}(X,Y)=\mathscr L_{M,U}(X,Y)$ and $\mathscr
K_{\gamma,U}(X,Y)=\mathscr K_{M,U}(X,Y)$. Moreover, if $X$ is normed
and $U=\big\{x\in X:\|x\|\leq\varepsilon\big\}$ then $\mathscr
L_{\gamma,\varepsilon}(X,Y)=\mathscr L_{\gamma,U}(X,Y)$ and
$\mathscr K_{\gamma,\varepsilon}(X,Y)=\mathscr K_{\gamma,U}(X,Y)$.
Thus, both $\mathscr L_{M,\varepsilon}(\mathbb R,\mathbb R)$ and
$\mathscr K_{M,\varepsilon}(\mathbb R,\mathbb R)$ are families of
demi-linear functions in $\mathbb R^\mathbb R$.

\begin{Definition}
([1,\ Def. 1.1]) A function $f:C_0^\infty(\Omega)\rightarrow\mathbb
C$ is called a demi-distribution if $f$ is continuous and
$f\in\mathscr L_{\gamma,U}(C_0^\infty(\Omega),\mathbb C)$ for some
$\gamma\in C(0)$ and $U\in\mathcal N(C_0^\infty(\Omega))$.
\end{Definition}

Let $C_0^\infty(\Omega)^{(\gamma,U)}$ (resp.,
$C_0^\infty(\Omega)^{[\gamma,U]}$) be the family of
demi-distributions which are functionals in $\mathscr
L_{\gamma,U}(C_0^\infty(\Omega),\mathbb C)$ (resp., $\mathscr
K_{\gamma,U}(C_0^\infty(\Omega),\mathbb C)$).

Let $C_0^\infty(\Omega)'$ be the family of usual distributions. Then
$$\mathscr D'(\Omega)=C_0^\infty(\Omega)'\subset C_0^\infty(\Omega)^{[\gamma,U]}
\subset C_0^\infty(\Omega)^{(\gamma,U)},\ \ \forall\,\gamma\in
C(0),\ U\in\mathcal N(C_0^\infty(\Omega))$$ and, in general,
$C_0^\infty(\Omega)^{[\gamma,U]}\backslash C_0^\infty(\Omega)'$
includes nonlinear functionals as many as usual distributions, at
least (see [1-3]).

Notice that the notations $\mathscr L_{\gamma,U}(X,Y)$, $\mathscr
K_{\gamma,U}(X,Y)$, $C_0^\infty(\Omega)^{(\gamma,U)}$ and
$C_0^\infty(\Omega)^{[\gamma,U]}$ always mean that $\gamma\in C(0)$
and $U\in\mathcal N(X)$ (resp., $\mathcal N(C_0^\infty(\Omega))$),
automatically. We also have similar understanding for $\mathscr
L_{M,U}(X,Y)$, $\mathscr K_{M,\varepsilon}(\mathbb C,\mathbb C)$,
etc.

\section{Continuity of Demi-distributions}
Throughout this paper, $n\in\mathbb N$ and $\Omega$ is a nonempty
open set in $\mathbb R^n$.

\begin{Definition}
$\mathcal{S}\subset\mathbb C^\Omega$, $\mathcal S\neq\emptyset$. For
$\xi\in\mathcal S$ and $f:\mathcal S\rightarrow\mathbb C$, let
$$supp\,\xi=\overline{\big\{x\in\Omega:\xi(x)\neq0\big\}}\cap\Omega,$$
$$supp\,f=\big\{x\in\Omega:\forall\,\mbox{open }G\subset\Omega\mbox{ with }x\in G
\,\exists\,\xi\in\mathcal S\mbox{ with }supp\,\xi\subset G\mbox{
such that }f(\xi)\neq0\big\}.$$
\end{Definition}

\begin{Lemma}
$\mathcal{S}\subset\mathbb C^\Omega$, $\mathcal S\neq\emptyset$. For
$\xi\in\mathcal S$ and $f\in\mathbb C^{\mathcal S}$, both
$supp\,\xi$ and $supp\,f$ are closed in $\Omega$ and so both
$\Omega\backslash supp\,\xi$ and $\Omega\backslash supp\,f$ are open
in $\mathbb R^n$.

Proof. Let $x_k\in supp\,f$ and $x_k\rightarrow x\in\Omega$. If
$x\not\in supp\,f$ then there is an open $G\subset\Omega$ with $x\in
G$ such that $f(\xi)=0$ for every $\xi\in\mathcal S$ with
$supp\,\xi\subset G$. But $x_k\in G$ eventually and so $x_k\not\in
supp\,f$ eventually. This contradiction shows that $x\in supp\,f$.\
\ $\square$
\end{Lemma}

\begin{Lemma}
Let $\xi\in C_0^\infty(\Omega)$ and $f\in\mathscr
L_{\gamma,U}(C_0^\infty(\Omega),\mathbb C)$. If $supp\,\xi\cap
supp\,f=\emptyset$, then $f(\xi)=0$.

Proof. If $\xi=0$ then $f(\xi)=f(0)=0$ by Def. 1.1.

Let $\xi\neq0$. Since $supp\,\xi\neq\emptyset$ and $supp\,f\cap
supp\,\xi=\emptyset$, $supp\,f\subsetneqq\Omega$. By Lemma 2.1, for
every $x\in supp\,\xi$ there is an open $G_x\subset\Omega\backslash
supp\,f$ such that $x\in G_x$ and $f(\eta)=0$ for all $\eta\in
C_0^\infty(G_x)$. Since $supp\,\xi$ is compact, there exist
$x_1,\cdots,x_m\in supp\,\xi$ such that
$supp\,\xi\subset\bigcup_{j=1}^m G_{x_j}$ and so $\xi\in
C_0^\infty(\bigcup_{j=1}^m G_{x_j})$. By Th. 1.4.4 of [4],
$\xi=\sum_{j=1}^m\xi_j$ where $\xi_j\in C_0^\infty(G_{x_j})$,
$j=1,2,\cdots,m$.

Pick a $p\in\mathbb N$ for which $\frac{1}{p}\xi_j\in U$,
$j=1,2,\cdots,m$. But each $\frac{1}{p}\xi_j\in C_0^\infty(G_{x_j})$
so $f(\frac{1}{p}\xi_j)=0$, $j=1,2,\cdots,m$. Therefore,
\begin{align*}
f(\xi)=f\Big(\sum_{j=1}^m\xi_j\Big)&=f\Big(\sum_{j=1}^{m-1}\xi_j+(p-1)\frac{1}{p}\xi_m+\frac{1}{p}\xi_m\Big)\\
&=r_1f\Big(\sum_{j=1}^{m-1}\xi_j+(p-1)\frac{1}{p}\xi_m\Big)+s_1f\Big(\frac{1}{p}\xi_m\Big)\\
&=r_1f\Big(\sum_{j=1}^{m-1}\xi_j+(p-2)\frac{1}{p}\xi_m+\frac{1}{p}\xi_m\Big)\\
&=r_1r_2\cdots
r_pf\Big(\sum_{j=1}^{m-1}\xi_j\Big)=\cdots=r_1r_2\cdots
r_{mp-1}f\Big(\frac{1}{p}\xi_1\Big)=0.\ \ \square
\end{align*}
\end{Lemma}

A linear functional $f:C_0^\infty(\Omega)\rightarrow\mathbb C$ is
continuous if and only if the condition (1.1) holds for $f$ [4, Th.
2.1.4]. However, for demi-linear functionals in $\mathscr
L_{\gamma,U}(C_0^\infty(\Omega),\mathbb C)$ the relation between
continuity and the condition (1.1) is quite complicated. First, we
show that many demi-distributions satisfy the condition (1.1).

\begin{Example}
(1) Let $n=1$ and $f(\xi)=\int_{-1}^1|\sin|\xi(x)||\,dx,\ \xi\in
C_0^\infty(\mathbb R)$. It is easy to see that $f$ is not linear but
$f\in C_0^\infty(\mathbb R)^{[\gamma,U]}$ where
$\gamma(t)=\frac{\pi}{2}t$ for $t\in\mathbb C$ and $U=\big\{\xi \in
C_0^\infty(\mathbb R):\sup_{|x|\leq1}|\xi(x)|\leq1\big\}$. For every
compact $K\subset\mathbb R$ and $\xi\in C_0^\infty(K)$,
$$\big|f(\xi)\big|=\int_{-1}^1\big|\sin|\xi(x)|\big|\,dx\leq\int_{-1}^1\big|\xi(x)\big|\,dx
\leq2\sup_{x\in
K}\big|\xi(x)\big|=2\sup_{K}\big|\partial^0\xi\big|.$$ Thus, $f$ is
a demi-distribution of order $0$. Moreover, the constant $C=2$ in
(1.1) is available for all compact $K\subset\mathbb R$.

(2) Pick a $f\in L_{loc}^1(\mathbb R^n)$ with $\sup_{x\in\mathbb
R^n}|f(x)|\leq M<+\infty$ and let
$$\big[f\big](\xi)=\int_{\mathbb R^n}\big|f(x)\xi(x)\big|\,dx,\ \ \xi\in C_0^\infty(\mathbb R^n)=\mathscr D.$$
Then $[f]$ is not linear but $[f]\in\mathscr D^{[\gamma,\mathscr
D]}$ for every $\gamma\in C(0)$ [1, Exam. 1.1(1)].

Let $K$ be a compact set in $\mathbb R^n$. Pick a cube $L\supset K$
for which $|L|=\int_L1\,dx<+\infty$. Then
$$\big|\big[f\big](\xi)\big|=\int_{\mathbb R^n}\big|f(x)\xi(x)\big|\,dx\leq M\int_L\big|\xi(x)\big|\,dx
\leq M|L|\sup_{K}\big|\partial^0\xi\big|,\ \ \forall\,\xi\in
C_0^\infty(K).$$ Thus, the condition (1.1) holds for the
demi-distribution $[f]$. If $f(x)=|x|=\sqrt{x_1^2+\cdots+x_n^2}$ for
all $x=(x_1,\cdots,x_n)\in\mathbb R^n$, then $supp\,[f]=\mathbb R^n$
is not compact but the condition (1.1) holds for $[f]$ and $[f]$ is
of order $0$.
\end{Example}

For $\gamma(t)=Mt$ where $M\geq1$ we will show that if $f\in
C_0^\infty(\Omega)^{[\gamma,U]}$ and $supp\,f$ is compact, then the
condition (1.1) holds for $f$ and,  in fact, $f$ has a more strong
property (see \S\,4, Th. 4.2), and Exam. 2.1(2) is interesting
because this example shows that the condition (1.1) can not imply
compactness of support. Moreover, the condition (1.1) fails to hold
for some $f\in C_0^\infty(\Omega)^{(\gamma,U)}\backslash
C_0^\infty(\Omega)^{[\gamma,U]}$, though $supp\,f$ is compact. We
show that (1.1) can be false even if $supp\,f=\{x_0\}$ is a
singleton.

\begin{Example}
For $a>0$ let $H_a(x)=1/a$ when $0<x<a$ and $H_a(x)=0$ otherwise.
Let $1\geq a_0>a_1>a_2>\cdots$ be a positive sequence with
$\sum_{j=0}^\infty a_j=a<+\infty$ and
$u=\lim_k(H_{a_0}\ast\cdots\ast H_{a_k})$. By Th. 1.3.5 of [4],
$u\in C_0^\infty(\mathbb R)$, $supp\,u\subset[0,a]$, $\int u\,dx=1$
and
$$\big|u^{(k)}(x)\big|\leq\frac{1}{2}\int\big|u^{(k+1)}(x)\big|\,dx\leq2^k/(a_0\cdots a_k),\ \ x\in\mathbb R,\ k=0,1,2\cdots.$$

Pick an $x_0\in(0,a)$ for which $u(x_0)=\sup_{x\in\mathbb R}u(x)>0$
and define a continuous $f:C_0^\infty(\mathbb R)\rightarrow\mathbb
R$ by $f(\xi)=e^{|\xi(x_0)|}-1$, $\xi\in C_0^\infty(\mathbb R)$.
Letting $\gamma(t)=et$ for $t\in\mathbb C$ and $U=\big\{\xi\in
C_0^\infty(\mathbb R):\sup_{0\leq x\leq a}|\xi(x)|\leq1\big\}$, it
is easy to see that $f\in\mathscr L_{\gamma,U}(C_0^\infty(\mathbb
R),\mathbb R)$ and so $f$ is a demi-distribution in
$C_0^\infty(\mathbb R)^{(\gamma,U)}$. Clearly, $supp\,f$ is compact
and, in fact, $supp\,f=\{x_0\}$.

Observe that $u\in C_0^\infty(\mathbb R)$. Then $mu\in
C_0^\infty(\mathbb R)$ and $supp\,(mu)=supp\,u\subset[0,a]$ for all
$m\in\mathbb N$, and
$$f(mu)=e^{|mu(x_0)|}-1=e^{mu(x_0)}-1=e^{c_m}mu(x_0)\mbox{ where }\lim_me^{c_m}=+\infty.$$
Now let $C>0$ and $k\in\mathbb N\cup\{0\}$. There is an
$m_0\in\mathbb N$ such that
$$e^{c_m}>C(k+1)\frac{2^k}{u(x_0)a_0a_1\cdots a_k},\ \forall\,m\geq m_0.$$
Then $|f(mu)|=e^{c_m}mu(x_0)>C(k+1)m\frac{2^k}{a_0a_1\cdots
a_k}>C\sum_{j=0}^km\frac{2^j}{a_0a_1\cdots a_j}\geq
C\sum_{j=0}^k\\|(mu)^{(j)}(x)|,\ \ \forall\,m\geq m_0,\ x\in\mathbb
R.$

Thus, for every $C>0$ and $k\in\mathbb N\cup\{0\}$ there exists
$m_0\in\mathbb N$ such that $mu\in C_0^\infty([0,a])$ for all $m\geq
m_0$ and
$$\big|f(mu)\big|>C\sum_{j=0}^k\sup_{x\in[0,a]}\big|(mu)^{(j)}(x)\big|,\ \ \forall\,m\geq m_0,$$
that is, the condition (1.1) fails to hold for $f$.
\end{Example}

However, for demi-distributions there is a simple condition implies
(1.1).

\begin{Theorem}
Let $f\in C_0^\infty(\Omega)^{(\gamma,U)}$. If there is an
$\varepsilon>0$ such that
$$\varepsilon\big|tf(\xi)\big|\leq\big|f(t\xi)\big|,\ \ \forall\,t>0,\ \xi\in C_0^\infty(\Omega),\leqno(2.1)$$
then the condition (1.1) holds for $f$.

Proof. If the conclusion fails, there is a compact $K\subset\Omega$
such that
$$\forall\,j\in\mathbb N\,\exists\,\xi_j\in C_0^\infty(K)
\mbox{ such that }\big|f(\xi_j)\big|>j\sum_{|\alpha|\leq
j}\sup_K\big|\partial^\alpha\xi_j\big|.\leqno(2.2)$$ Then
$|f(\xi_j)|>0$, $\xi_j\neq0$, $\xi_j(x_j)\neq0$ for some $x_j\in K$
and
$$\sum_{|\alpha|\leq j}\sup_K\big|\partial^\alpha\xi_j\big|\geq\sup_K\big|\partial^0\xi_j\big|
\geq\big|\xi_j(x_j)\big|>0,\ \ j=1,2,3,\cdots.$$ It follows from
(2.1) and (2.2) that
$$\varepsilon<\varepsilon\Big|\frac{f(\xi_j)}{j\sum_{|\alpha|\leq j}\sup_K|\partial^\alpha\xi_j|}\Big|\leq
\Big|f(\frac{\xi_j}{j\sum_{|\alpha|\leq
j}\sup_K|\partial^\alpha\xi_j|})\Big|,\ j=1,2,3,\cdots.\leqno(2.3)$$

Let $\beta$ be a multi-index. Then
$$\sup_{x\in K}\Big|\partial^\beta\big(\frac{\xi_j}{j\sum_{|\alpha|\leq j}
\sup_K|\partial^\alpha\xi_j|}\big)(x)\Big|\leq\frac{1}{j},\ \
\forall\,j\geq|\beta|,$$
$$\lim_{j\rightarrow\infty}\sup_K\Big|\partial^\beta\big(\frac{\xi_j}{j\sum_{|\alpha|\leq j}
\sup_K|\partial^\alpha\xi_j|}\big)\Big|=0,\ \
\forall\,\mbox{multi-index}\ \beta.\leqno(2.4)$$ Observing
$\big\{\xi_j/j\sum_{|\alpha|\leq
j}\sup_K|\partial^\alpha\xi_j|\big\}\subset C_0^\infty(K)$, it
follows from (2.4) that
$$\frac{\xi_j}{j\sum_{|\alpha|\leq j}\sup_K|\partial^\alpha\xi_j|}\rightarrow0\mbox{ in }C_0^\infty(\Omega).$$
Since $f\in C_0^\infty(\Omega)^{(\gamma,U)}$ is continuous,
$f\big(\frac{\xi_j}{j\sum_{|\alpha|\leq
j}\sup_K|\partial^\alpha\xi_j|}\big)\rightarrow0$ but this
contradicts (2.3).\ \ $\square$
\end{Theorem}

\begin{Theorem}
Every nonzero usual distribution $f\in\mathscr D'(\Omega)$ produces
uncountably many of nonlinear demi-distributions satisfying the
conditions (2.1) and (1.1).

Proof. Let $f\in\mathscr D'(\Omega)$, $f\neq0$. There is
$U\in\mathcal N(C_0^\infty(\Omega))$ such that $|f(\eta)|\leq1$ for
all $\eta\in U$.

Pick a nonlinear continuous $h:\mathbb R\rightarrow\mathbb R$ such
that
$$h(x)=x\mbox{ when }|x|\leq1,\ \ \frac{1}{2}\leq\frac{h(b)-h(a)}{b-a}\leq1\mbox{ when }a<b.$$
Clearly, $\mathbb R^\mathbb R$ includes uncountably many of this
kind functions. For every $\varepsilon\in(0,\frac{1}{2}]$ and
$x\in\mathbb R$ we have
$$\varepsilon|x|\leq\big|h(x)\big|\leq|x|,\ \ \varepsilon|x|\leq\big|h(|x|)\big|=h\big(|x|\big)\leq|x|.$$

By Th. 1.1 of [2], $h\in\mathscr K_{1,1}(\mathbb R,\mathbb R)$, that
is, for $x\in\mathbb R$ and $u,t\in[-1,1]$ we have that
$h(x+tu)=h(x)+sh(u)$ where $|s|\leq|t|$. Then for $\xi\in
C_0^\infty(\Omega)$, $\eta\in U$ and $|t|\leq1$,
$h(|f(\xi+t\eta)|)=h(|f(\xi)+tf(\eta)|)=h(|f(\xi)|+s|f(\eta)|)$
where $|s|\leq|t|\leq1$ and, therefore,
$$h\big(|f(\xi+t\eta)|\big)=h\big(|f(\xi)|\big)+s'h\big(|f(\eta)|\big),\ \ |s'|\leq|s|\leq|t|\leq|\gamma(t)|,
\ \forall\,\gamma\in C(0).$$ This shows that $h(|f(\cdot)|)\in
C_0^\infty(\Omega)^{[\gamma,U]}\backslash\mathscr D'(\Omega)$,
$\forall\,\gamma\in C(0)$.

Let $0<\varepsilon\leq\frac{1}{2}$. Then
$$\varepsilon\big|th\big(|f(\xi)|\big)\big|\leq\varepsilon\big|tf(\xi)\big|=\varepsilon\big|f(t\xi)\big|
\leq\big|h\big(|f(t\xi)|\big)\big|,\ \ \forall\,t\in\mathbb R,\
\xi\in C_0^\infty(\Omega).$$ So (2.1) holds for $h(|f(\cdot)|)$. By
Th. 2.1, $h(|f(\cdot)|)$ satisfies the condition (1.1), and the
usual distribution $f$ produces uncountably many of this kind
nonlinear demi-distributions.\ \ $\square$
\end{Theorem}

We also are interested in the converse implications.

\begin{Theorem}
Let $f\in\mathscr L_{\gamma,U}(C_0^\infty(\Omega),\mathbb C)$. If
the condition (1.1) holds for $f$, then $f$ is sequentially
continuous.

Proof. Suppose $\xi_j\rightarrow\xi$ in $C_0^\infty(\Omega)$. Then
$\xi_j-\xi\rightarrow0$ and there is a compact $K\subset\Omega$ such
that $supp\,(\xi_j-\xi)\subset K$ for all $j$ [4, p.35]. Moreover,
there exist sequences $t_j\rightarrow0$ in $\mathbb C$ and
$\eta_j\rightarrow0$ in $C_0^\infty(\Omega)$ such that
$\xi_j-\xi=t_j\eta_j$ for all $j$ [6, Exam. 2]. We may assume that
$|t_j|\leq1$ and $\eta_j\in U$ for all $j$. Then
$$f(\xi_j)-f(\xi)=f(\xi+\xi_j-\xi)-f(\xi)=f(\xi+t_j\eta_j)-f(\xi)=(r_j-1)f(\xi)+s_jf(\eta_j),$$
where $|r_j-1|\leq|\gamma(t_j)|\rightarrow0$ and
$|s_j|\leq|\gamma(t_j)|\rightarrow0$.

If $t_j=0$ then $f(\xi_j)=f(\xi+t_j\eta_j)=f(\xi)$ so we may assume
that $t_j\neq0$ for all $j$. Then
$supp\,\eta_j=supp\,(t_j\eta_j)=supp\,(\xi_j-\xi)\subset K$ for all
$j$ and by the condition (1.1) there exist $C>0$ and $k\in\mathbb
N\cup\{0\}$ such that
$$\big|f(\eta_j)\big|\leq C\sum_{|\alpha|\leq k}\sup_K\big|\partial^\alpha\eta_j\big|\rightarrow0
\mbox{ as }j\rightarrow\infty\mbox{ since }\eta_j\rightarrow0\mbox{
in }C_0^\infty(\Omega),$$ so $f(\eta_j)\rightarrow0$. Thus,
$f(\xi_j)-f(\xi)=(r_j-1)f(\xi)+s_jf(\eta_j)\rightarrow0$. \ \
$\square$
\end{Theorem}

Let $\gamma_0\in C(0)$, $\gamma_0(t)=t$ for $t\in\mathbb C$. For
every $U\in\mathcal N(C_0^\infty(\Omega))$ the family $\mathscr
K_{1,U}(C_0^\infty(\Omega),\mathbb C)=\mathscr
K_{\gamma_0,U}(C_0^\infty(\Omega),\mathbb C)$ includes all linear
functionals and much more nonlinear functionals, e.g., for every
nonzero linear $f:C_0^\infty(\Omega)\rightarrow\mathbb C$,
$|f(\cdot)|$ is nonlinear but $|f(\cdot)|\in\mathscr
K_{1,U}(C_0^\infty(\Omega),\mathbb C),\ \forall\,U\in\mathcal
N(C_0^\infty(\Omega))$.

\begin{Theorem}
If $f\in\mathscr K_{1,U}(C_0^\infty(\Omega),\mathbb C)$ and the
condition (1.1) holds for $f$, then $f$ is continuous and so $f$ is
a demi-distribution in $C_0^\infty(\Omega)^{[\gamma_0,U]}$.

Proof. By Th. 2.3, $f$ is sequentially continuous. Since
$C_0^\infty(\Omega)$ is bornological, $C_0^\infty(\Omega)$ is
C-sequential and $f$ is continuous by Th. 1.1 of [1]. \ \ $\square$
\end{Theorem}

In \S\,4 we will improve this result (see Cor. 4.2).

\section{Extensions of Demi-distributions}
For every $C\geq1$ and $\varepsilon>0$, $\mathscr
K_{C,\varepsilon}(\mathbb C,\mathbb C)$ includes uncountably many
nonlinear functionals, and $h\circ f$ is a demi-distribution in
$\Omega$ for every $h\in\mathscr L_{\gamma,\varepsilon}(\mathbb
C,\mathbb C)$ and $f\in\mathscr D'(\Omega)$ (see [1, Th. 1.5, Cor.
1.1]).

\begin{Theorem}
If $h\in\mathscr L_{\gamma,\varepsilon}(\mathbb C,\mathbb C)$,
$h\neq0$ and $f\in\mathscr D'(\Omega)$, then $h\circ f$ is a
demi-distribution in $\Omega$ and
$$supp\,\big(h\circ f\big)=supp\,f.$$

Proof. If $x\in\Omega\backslash supp\,f$ then there is an open
$N_x\subset\Omega$ such that $x\in N_x$ and $f(\eta)=0$ for all
$\eta\in C_0^\infty(N_x)$. Then $(h\circ f)(\eta)=h(f(\eta))=h(0)=0$
when $\eta\in C_0^\infty(N_x)$ and so $x\not\in supp\,(h\circ f)$.
Thus $supp\,(h\circ f)\subset supp\,f$.

If $u\in\mathbb C$ such that $0<|u|<\varepsilon$ and $h(u)=0$, then
for every $z\in\mathbb C$ there is a $p\in\mathbb N$ such that
$\frac{1}{p}|\frac{z}{u}|\leq1$ and
$h(z)=h(p\frac{z}{pu}u)=h((p-1)\frac{z}{pu}u+\frac{z}{pu}u)=r_1h((p-1)\frac{z}{pu}u)+s_1h(u)
=r_1h((p-1)\frac{z}{pu}u)=\cdots=r_1r_2\cdots r_{p-1}s_ph(u)=0$.
This contradicts that $h\neq0$. Hence $h(u)\neq0$ when
$0<|u|<\varepsilon$.

Let $x\in supp\,f$ and $N_x$ an open neighborhood of $x$ such that
$N_x\subset\Omega$. Then $f(\eta)\neq0$ for some $\eta\in
C_0^\infty(N_x)$ and $0<|\frac{1}{p}f(\eta)|<\varepsilon$ for some
$p\in\mathbb N$. Observing $f$ is a usual distribution,
$\frac{1}{p}\eta\in C_0^\infty(N_x)$ and $(h\circ
f)(\frac{1}{p}\eta)=h[f(\frac{1}{p}\eta)]=h[\frac{1}{p}f(\eta)]\neq0$.
This shows that $x\in supp\,(h\circ f)$ so $supp\,f\subset
supp\,(h\circ f)$.\ \ $\square$
\end{Theorem}

For $M\subset\Omega$ let $C_0(M)=\big\{\xi\in\mathbb
C^\Omega:\xi\mbox{ is continuous}, supp\,\xi\mbox{ is
compact}\big\}$. We have an analogue of Th. 1.4.4 of [4] as follows.

\begin{Lemma}
Let $\Omega_1,\cdots,\Omega_k$ be open sets in $\Omega$ and let
$\xi\in C_0(\bigcup_{1}^k\Omega_j)$. Then one can find $\xi_j\in
C_0(\Omega_j)$, $j=1,2,\cdots,k$, such that $\xi=\sum_1^k\xi_j$. If
$\xi\geq0$ one can take all $\xi_j\geq0$.

Proof. If $x\in supp\,\xi$ then $x\in\Omega_j$ for some
$j\in\{1,2,\cdots,k\}$ and there is a compact neighborhood of $x$
contained in $\Omega_j$. Since $supp\,\xi$ is compact, a finite
number of such neighborhoods can be chosen which cover all of
$supp\,\xi$. Hence $supp\,\xi\subset \bigcup_1^k K_j$ where each
$K_j$ is compact and $K_j\subset\Omega_j$.

By Th. 1.4.1 of [4], there is $\mathcal {X}_j\in
C_0^\infty(\Omega_j)$ such that $0\leq\mathcal X_j\leq1$ and
$\mathcal X_j=1$ in a neighborhood of $K_j$, $j=1,2,\cdots,k$. Let
$$\xi_1=\xi\mathcal X_1,\ \xi_2=\xi\mathcal X_2(1-\mathcal X_1),\
\cdots,\ \xi_k=\xi\mathcal X_k(1-\mathcal X_1)\cdots(1-\mathcal
X_{k-1}),$$ then each $supp\,\xi_j\subset supp\,\mathcal
X_j\subset\Omega_j$  and $\xi=\sum_1^k\xi_j$.\ \ $\square$
\end{Lemma}

Let $\mathcal S\subset C(\Omega)$, $M\subset\Omega$ and $\mathcal
S(M)=\big\{\xi\in\mathcal S:supp\,\xi\subset M\big\}$. For a
function $f:\mathcal S\rightarrow\mathbb C$ define $f_M:\mathcal
S(M)\rightarrow\mathbb C$ by $f_M(\xi)=f(\xi),\
\forall\,\xi\in\mathcal S(M)$, and $f_M$ is called the restriction
of $f$ to $M$.

We now improve Th. 2.2.1 of [4].

\begin{Theorem}
Let $f\in \mathscr L_{\gamma,U}(C_0(\Omega),\mathbb C)$. If every
point in $\Omega$ has a neighborhood to which the restriction of $f$
is $0$, then $f=0$. The same fact is valid for $f\in\mathscr
L_{\gamma,U}(C_0^\infty(\Omega),\mathbb C)$.

Proof. Let $\xi\in C_0(\Omega)$. Since $supp\,\xi$ is compact, there
exist $x_1,\,\cdots,\,x_m\in supp\,\xi$ such that
$supp\,\xi\subset\bigcup_1^mN_j$ where $N_j$ is an open neighborhood
of $x_j$ such that $f_{N_j}=0$, $1\leq j\leq m$. Then $\xi\in
C_0(\bigcup_1^mN_j)$ and $\xi=\sum_1^m\xi_j$ where each $\xi_j\in
C_0(N_j)$ by Lemma 3.1.

Pick a $p\in\mathbb N$ such that $\frac{1}{p}\xi$,
$\frac{1}{p}\xi_j\in U$, $j=1,2,\cdots,m$. Then
$\frac{1}{p}\xi=\sum_1^m\frac{1}{p}\xi_j$ and
$f(\frac{1}{p}\xi)=f(\sum_1^m\frac{1}{p}\xi_j)=r_1f(\sum_{1}^{m-1}\frac{1}{p}\xi_j)+s_1f(\frac{1}{p}\xi_m)
=r_1f(\sum_1^{m-1}\frac{1}{p}\xi_j)=\cdots=r_1r_2\cdots
r_{m-1}f(\frac{1}{p}\xi_1)=0$ since each $\frac{1}{p}\xi_j\in
C_0(N_j)$ and $f(\frac{1}{p}\xi_j)=f_{N_j}(\frac{1}{p}\xi_j)=0$.
Thus,
$f(\xi)=f(p\frac{1}{p}\xi)=t_1f((p-1)\frac{1}{p}\xi)=\cdots=t_1t_2\cdots
t_{p-1}f(\frac{1}{p}\xi)=0$.

The same conclusion can be obtained for $f\in\mathscr
L_{\gamma,U}(C_0^\infty(\Omega),\mathbb C)$ using Th. 1.4.4 of [4]
instead of Lemma 3.1.\ \ $\square$
\end{Theorem}

\begin{Theorem}
If $f\in\mathscr L_{\gamma,U}(C_0(\Omega),\mathbb C)$ (resp.,
$\mathscr L_{\gamma,U}(C_0^\infty(\Omega),\mathbb C)$) and $\xi\in
C_0(\Omega)$ (resp., $C_0^\infty(\Omega)$) such that $supp\,f\cap
supp\,\,\xi=\emptyset$, then $f(\xi)=0$.

Proof. Let $x\in supp\,\xi$. Since $x\not\in supp\,f$ and $supp\,f$
is closed in $\Omega$ by Lemma 2.1, there is an open neighborhood
$N_x$ of $x$ such that $N_x\subset\Omega\backslash supp\,f$ and the
restriction $f_{N_x}=0$. Since $supp\,\xi$ is compact, there exist
$x_1,\cdots,x_m\in supp\,\xi$ such that
$supp\,\xi\subset\bigcup_1^mN_{x_j}$ and $f_{N_{x_j}}=0$,
$j=1,2,\cdots,m$. Then $\xi\in C_0(\bigcup_1^mN_{x_j})$ (resp.,
$C_0^\infty(\bigcup_1^m N_{x_j})$) and $\xi=\sum_1^m\xi_j$ where
$\xi_j\in C_0(N_{x_j})$ (resp., $C_0^\infty(N_{x_j})$) by Lemma 3.1
(resp., Th. 1.4.4 of [4]), $j=1,2,\cdots,m$.

Now $f(\xi)=0$ as in the proof of Th. 3.2.\ \ $\square$
\end{Theorem}

Note that Th. 3.3 is not a consequence of Th. 3.2 because for
$f\neq0$ and $x\in supp\,f$ the restriction $f_{N_x}\neq0$ when
$N_x$ is a neighborhood of $x$.

\begin{Definition}
Let $f\in\mathscr K_{\gamma,U}(C_0^\infty(\Omega),\mathbb C)$ and
$\xi\in C^\infty(\Omega)$. We say that $\xi=\xi_0+\xi_1$ is a
$f$-decomposition of $\xi$ if $\xi_0\in C_0^\infty(\Omega)$ and
$supp\,\xi_1\cap supp\,f=\emptyset$.
\end{Definition}

Observe that for every compact $K\subset\Omega$ there is a $\mathcal
X\in C_0^\infty(\Omega)$ such that $0\leq\mathcal X\leq1$ and
$\mathcal X=1$ in a neighborhood of $K$.

\begin{Lemma}
Let $f\in\mathscr K_{\gamma,U}(C_0^\infty(\Omega),\mathbb C)$ and
$\xi\in C^\infty(\Omega)$ such that $supp\,\xi\cap supp\,f$ is
compact. If $K$ is a compact subset of $\Omega$ such that
$supp\,\xi\cap supp\,f\subseteq K$ and $\mathcal X\in
C_0^\infty(\Omega)$ for which $\mathcal X=1$ in a neighborhood of
$K$, then $\xi=\mathcal X\xi+(1-\mathcal X)\xi$ is a
$f$-decomposition of $\xi$: $\mathcal X\xi\in C_0^\infty(\Omega)$,
$supp\,[(1-\mathcal X)\xi]\cap supp\,f=\emptyset$.

Proof. Since $supp\,(\mathcal X\xi)\subset supp\,\mathcal X$,
$\mathcal X\xi\in C_0^\infty(\Omega)$. There is an open
$G\subset\Omega$ such that $K\subset G$ and $\mathcal X=1$ in $G$.
If $[(1-\mathcal X)\xi](x)=(1-\mathcal X(x))\xi(x)\neq0$, then $x\in
supp\,\xi\cap(\Omega\backslash G)$ and so $supp\,[(1-\mathcal
X)\xi]\subset supp\,\xi\cap(\Omega\backslash G)\subset
supp\,\xi\cap(\Omega\backslash K)\subset
supp\,\xi\cap[\Omega\backslash(supp\,\xi\cap
supp\,f)]\subset\Omega\backslash supp\,f$, i.e., $supp\,[(1-\mathcal
X)\xi]\cap supp\, f=\emptyset$.\ \ $\square$
\end{Lemma}

\begin{Corollary}
Let $f\in\mathscr K_{\gamma,U}(C_0^\infty(\Omega),\mathbb C)$ and
$\xi\in C_0^\infty(\Omega)$. If $K$ is a compact subset of $\Omega$
such that $supp\,\xi\cap supp\,f\subseteq K$ and $\mathcal X\in
C_0^\infty(\Omega)$ for which $\mathcal X=1$ in a neighborhood of
$K$, then $f(\xi)=f(\mathcal X\xi)$.

Proof. By Lemma 3.2, $\xi=\mathcal X\xi+(1-\mathcal X)\xi$ is a
$f$-decomposition of $\xi$ and $supp\,[(1-\mathcal X)\xi]\cap
supp\,f=\emptyset$. But $\xi\in C_0^\infty(\Omega)$ so $(1-\mathcal
X)\xi\in C_0^\infty(\Omega)$.

Pick a $p\in\mathbb N$ such that $\frac{1}{p}(1-\mathcal X)\xi\in
U$. Then $supp\,[\frac{1}{p}(1-\mathcal X)\xi]\cap
supp\,f=supp\,[(1-\mathcal X)\xi]\cap supp\,f=\emptyset$ and
$f[\frac{1}{p}(1-\mathcal X)\xi]=0$ by Lemma 2.2 (see also Th. 3.3).
Hence

$f\big(\xi\big)=f\big[\mathcal X\xi+p\frac{1}{p}(1-\mathcal
X)\xi\big]=f\big[\mathcal X\xi+(p-1)\frac{1}{p}(1-\mathcal
X)\xi\big]=\cdots=f\big[\mathcal X\xi+\frac{1}{p}(1-\mathcal
X)\xi\big]=f\big(\mathcal X\xi\big)+sf\big[\frac{1}{p}(1-\mathcal
X)\xi\big]=f\big(\mathcal X\xi\big)$.\ \ $\square$
\end{Corollary}

\begin{Theorem}
Let $f\in\mathscr K_{\gamma,U}(C_0^\infty(\Omega),\mathbb C)$ and
$\xi\in C^\infty(\Omega)$. If both $\xi=\xi_0+\xi_1$ and
$\xi=\eta_0+\eta_1$ are $f$-decompositions of $\xi$, where $\xi_0$,
$\eta_0\in C_0^\infty(\Omega)$ and $supp\,\xi_1\cap
supp\,f=supp\,\eta_1\cap supp\,f=\emptyset$, then
$f(\xi_0)=f(\eta_0)$.

Proof. Let $\mathcal X=\xi_0-\eta_0$. Then $supp\,\mathcal X\subset
supp\,\xi_0\cup supp\,\eta_0$ and so $\mathcal X\in
C_0^\infty(\Omega)$. Since $\eta_1-\xi_1=\xi_0-\eta_0=\mathcal X$ so
$supp\,\mathcal X\subset supp\,\xi_1\cup supp\,\eta_1$,
$supp\,\mathcal X\cap supp\,f\subset(supp\,\xi_1\cap
supp\,f)\bigcup(supp\,\eta_1\cap supp\,f)=\emptyset$.

Pick a $p\in\mathbb N$ for which $\frac{1}{p}\mathcal X\in U$. Then
$supp\,(\frac{1}{p}\mathcal X)\cap(supp\,f)=\emptyset$ and so
$f(\frac{1}{p}\mathcal X)=0$ by Lemma 2.2. Then
$\xi_0=\eta_0+\mathcal X$ and
$$f\big(\xi_0\big)=f\big(\eta_0+p\frac{1}{p}\mathcal X\big)=
f\big(\eta_0+(p-1)\frac{1}{p}\mathcal
X\big)=\cdots=f\big(\eta_0+\frac{1}{p}\mathcal
X\big)=f\big(\eta_0\big).\ \ \square$$
\end{Theorem}

\begin{Corollary}
Let $f\in\mathscr K_{\gamma,U}(C_0^\infty(\Omega),\mathbb C)$ and
$\xi,\eta\in C_0^\infty(\Omega)$. If $supp\,(\xi-\eta)\cap
supp\,f=\emptyset$, then $f(\xi)=f(\eta)$.

Proof. Both $\xi=\xi+0$ and $\xi=\eta+(\xi-\eta)$ are
$f$-decompositions of $\xi$ so $f(\xi)=f(\eta)$.\ \ $\square$
\end{Corollary}

By Lemma 3.2 and Th. 3.4 we have

\begin{Definition}
For $f\in\mathscr K_{\gamma,U}(C_0^\infty(\Omega),\mathbb C)$ let
$$\mathcal S(f)=\big\{\xi\in C^\infty(\Omega):supp\,\xi\cap supp\,f\mbox{ is compact}\big\},$$
and define $\widetilde{f}:\mathcal S(f)\rightarrow\mathbb C$ by
$$\widetilde{f}(\xi)=f(\xi_0)\mbox{ when }\xi
=\xi_0+\xi_1\mbox{ is a }f\mbox{-decomposition of }\xi\in\mathcal
S(f).$$
\end{Definition}

We say that $\widetilde{f}$ is the canonical extension of $f$. If
$supp\,f$ is compact then $\mathcal S(f)=C^\infty(\Omega)$ and
$\widetilde{f}$ is defined on $C^\infty(\Omega)$.

\begin{Theorem}
Let $f\in\mathscr K_{\gamma,U}(C_0^\infty(\Omega),\mathbb C)$. Then
$\mathcal S(f)$ is a vector subspace of $C^\infty(\Omega)$ and
$C_0^\infty(\Omega)\subset\mathcal S(f)$. Moreover,
$$\widetilde{f}(\xi)=f(\xi),\ \xi\in C_0^\infty(\Omega),$$
$$\widetilde{f}(\xi)=0\mbox{ when }\xi\in C^\infty(\Omega)\mbox{ but }supp\,\xi\cap supp\,f=\emptyset.$$

Proof. If $\xi,\eta\in\mathcal S(f)$ and $t\in\mathbb C$, then
$(supp\,\xi\cap supp\,f)\bigcup(supp\,\eta\cap supp\,f)$ is compact
and $supp\,(\xi+t\eta)\cap supp\,f\subset(supp\,\xi\cup
supp\,\eta)\bigcap supp\,f$. This shows that $\xi+t\eta\in\mathcal
S(f)$. If $\xi\in C_0^\infty(\Omega)$ then $supp\,\xi\cap
supp\,f\subset supp\,\xi$ so $\xi\in\mathcal S(f)$, and
$\widetilde{f}=f(\xi)$ since $\xi=\xi+0$ is a $f$-decomposition of
$\xi$.

If $\xi\in C^\infty(\Omega)$ but $supp\,\xi\cap supp\, f=\emptyset$,
then $\xi=0+\xi$ is a $f$-decomposition of $\xi$ and
$\widetilde{f}(\xi)=f(0)=0$ by Def. 1.1.\ \ $\square$
\end{Theorem}

Recall that $C^\infty(\Omega)$ is a Fr\'echet space.

\begin{Lemma}
Let $\eta\in C^\infty(\Omega)$ and $T_\eta(\xi)=\eta\xi$ for $\xi\in
C^\infty(\Omega)$. Then $T_\eta:C^\infty(\Omega)\rightarrow
C^\infty(\Omega)$ is a continuous linear operator.

Proof. Let $\xi_v\rightarrow0$ in $C^\infty(\Omega)$. For every
compact $K\subset\Omega$ and $k\in\mathbb N\cup\{0\}$, it follows
from the Leibniz formula that
\begin{align*}
\big\|\eta\xi_v\big\|_{K,k}&=\sum_{|\alpha|\leq
k}\sup_{K}\big|\partial^\alpha(\eta\xi_v)\big|\\
&\leq
C\sum_{|\alpha|+|\beta|\leq k}\sup_K\big|\partial^\alpha\xi_v\big|\big|\partial^\beta\eta\big|\\
&\leq C\max_{|\beta|\leq
k}\sup_{K}\big|\partial^\beta\eta\big|\sum_{|\alpha|\leq
k}\sup_K\big|\partial^\alpha\xi_v\big|\rightarrow0.\ \ \square
\end{align*}
\end{Lemma}

The topology of the inductive limit $C_0^\infty(\Omega)$ is strictly
stronger than the topology of the subspace $C_0^\infty(\Omega)$ of
$C^\infty(\Omega)$. So the following fact is interesting and useful
for further discussions.

\begin{Lemma}
Let $\mathcal X\in C_0^\infty(\Omega)$ and $T_\mathcal
X(\xi)=\mathcal X\xi$ for $\xi\in C^\infty(\Omega)$. Then
$T_\mathcal X:C^\infty(\Omega)\rightarrow C_0^\infty(\Omega)$ is a
continuous linear operator.

Proof. Let $\xi_v\rightarrow0$ in $C^\infty(\Omega)$. Since
$\mathcal X\in C_0^\infty(\Omega)$, $supp\,\mathcal X$ is compact
and $supp\,(\mathcal X\xi_v)\subset supp\,\mathcal X$ for all $v$.
By Lemma 3.3, $\mathcal X\xi_v\rightarrow0$ in $C^\infty(\Omega)$
and so for every compact $K\subset\Omega$ and every multi-index
$\alpha$, $\lim_v\sup_K|\partial^\alpha(\mathcal
X\xi_v)|\leq\lim_v\sum_{|\beta|\leq|\alpha|}\\\sup_K|\partial^\beta(\mathcal
X\xi_v)|=\lim_v\|\mathcal X\xi_v\|_{K,|\alpha|}=0$. Thus $\mathcal
X\xi_v\rightarrow0$ in $C_0^\infty(\Omega)$ and $T_\mathcal
X:C^\infty(\Omega)\rightarrow C_0^\infty(\Omega)$ is continuous
because $C^\infty(\Omega)$ is a Fr\'echet space.\ \ $\square$
\end{Lemma}

\begin{Lemma}
Let $X,\ Y$ be topological vector spaces and $f\in\mathscr
L_{\gamma,U}(X,Y)$. Then $f$ is continuous if and only if $f$ is
continuous at $0\in X$.

Proof. Suppose that $f$ is continuous at $0\in X$. Let $x\in X$ and
$V\in\mathcal N(Y)$. Pick a balanced $W\in\mathcal N(Y)$ for which
$W+W\subset V$.

There is a balanced $U_0\in\mathcal N(X)$ such that $U_0\subset U$
and $f(U_0)\subset W$. Since $\lim_{t\rightarrow0}\gamma(t)=0$,
there is a $p\in\mathbb N$ for which $|\gamma(\frac{1}{p})|<1$ and
$\gamma(\frac{1}{p})f(x)\in W$. If $z\in x+\frac{1}{p}U_0$, then
$p(z-x)\in U_0\subset U$ and
$f(z)-f(x)=f(x+z-x)-f(x)=f[x+\frac{1}{p}p(z-x)]-f(x)=rf(x)+sf[p(z-x)]-f(x)=(r-1)f(x)+sf[p(z-x)]$,
where $|r-1|\leq|\gamma(\frac{1}{p})|<1$ and
$|s|\leq|\gamma(\frac{1}{p})|<1$.

If $\gamma(\frac{1}{p})=0$ then $r-1=s=0$ so $f(z)-f(x)=0\in V$. If
$\gamma(\frac{1}{p})\neq0$, then
$(r-1)f(x)=\frac{r-1}{\gamma(1/p)}\gamma(1/p)f(x)\in
\frac{r-1}{\gamma(1/p)}W\subset W$ and $sf[p(z-x)]\in sf(U_0)\subset
sW\subset W$. So $f(z)-f(x)\in W+W\subset V$. Thus,
$\frac{1}{p}U_0\in\mathcal N(X)$ and $f(x+\frac{1}{p}U_0)\subset
f(x)+V$, i.e., $f$ is continuous at $x$.\ \ $\square$
\end{Lemma}

Recall that $C_0^\infty(\Omega)^{[\gamma,U]}=\big\{f\in\mathscr
K_{\gamma,U}(C_0^\infty(\Omega),\mathbb C):f\mbox{ is
continuous}\big\}$ is the family of demi-distributions, and
$C_0^\infty(\Omega)^{[\gamma,U]}$ is a large extension of the family
$\mathscr D'(\Omega)(=C_0^\infty(\Omega)')$ of usual distributions.

\begin{Theorem}
Let $f\in C_0^\infty(\Omega)^{[\gamma,U]}$ such that $supp\,f$ is
compact. Then there is a $V\in\mathcal N(C^\infty(\Omega))$ such
that the canonical extension $\widetilde{f}\in
C^\infty(\Omega)^{[\gamma,V]}$ and $supp\,\widetilde{f}=supp\,f$.

Proof. Since $supp\,f$ is compact, $\mathcal S(f)=C^\infty(\Omega)$
and the canonical extension $\widetilde{f}$ of $f$ is defined on
$C^\infty(\Omega)$. Pick a $\mathcal X\in C_0^\infty(\Omega)$ such
that $\mathcal X=1$ in a neighborhood of $supp\,f$. By Lemma 3.2 and
Def. 3.2, $\widetilde{f}(\xi)=f(\mathcal X\xi)$, $\forall\,\xi\in
C^\infty(\Omega)$.

By Lemma 3.4, $V=\big\{\xi\in C^\infty(\Omega):\mathcal X\xi\in
U\big\}\in\mathcal N(C^\infty(\Omega))$. If $\xi\in
C^\infty(\Omega)$, $\eta\in V$ and $|t|\leq1$, then
$\widetilde{f}(\xi+t\eta)=f(\mathcal X\xi+t\mathcal
X\eta)=f(\mathcal X\xi)+sf(\mathcal
X\eta)=\widetilde{f}(\xi)+s\widetilde{f}(\eta)$ where
$|s|\leq|\gamma(t)|$. Thus $\widetilde{f}\in\mathscr
K_{\gamma,V}(C^\infty(\Omega),\mathbb C)$.

Let $\xi_v\rightarrow0$ in $C^\infty(\Omega)$. By Lemma 3.4,
$\mathcal X\xi_v\rightarrow0$ in $C_0^\infty(\Omega)$ and so
$\widetilde{f}(\xi_v)=f(\mathcal X\xi_v)\rightarrow
f(0)=0=\widetilde{f}(0)$. This shows that $\widetilde{f}$ is
continuous at $0\in C^\infty(\Omega)$ since $C^\infty(\Omega)$ is a
Fr\'echet space. Thus $\widetilde{f}$ is continuous by Lemma 3.5.

Let $x\in\Omega\backslash supp\,f$. There is an open
$N_x\subset\Omega\backslash supp\,f$ such that $x\in N_x$ and
$f(\eta)=0,\ \forall\,\eta\in C_0^\infty(N_x)$. If $\xi\in
C^\infty(N_x)$, then $supp\,(\mathcal X\xi)\subset supp\,\xi\subset
N_x$ so $\widetilde{f}(\xi)=f(\mathcal X\xi)=0$. Thus, $x\not\in
supp\,\widetilde{f}$ and so $supp\,\widetilde{f}\subset supp\,f$.
Conversely, if $x\in\Omega\backslash supp\,\widetilde{f}$ then there
is an open $N_x\subset supp\,\widetilde{f}$ such that
$\widetilde{f}(\xi)=0$ for all $\xi\in C^\infty(N_x)$ so
$f(\eta)=\widetilde{f}(\eta)=0,\ \forall\,\eta\in C_0^\infty(N_x)$.
Then $x\not\in supp\,f$ and so $supp\,f\subset
supp\,\widetilde{f}$.\ \ $\square$
\end{Theorem}

\begin{Theorem}
Let $f\in C^\infty(\Omega)^{[\gamma,V]}$ and define
$f_0:C_0^\infty(\Omega)\rightarrow\mathbb C$ by $f_0(\xi)=f(\xi)$
for $\xi\in C_0^\infty(\Omega)$. Then $U=\big\{\eta\in
C_0^\infty(\Omega):\eta\in V\big\}\in\mathcal N(C_0^\infty(\Omega))$
and $f_0\in C_0^\infty(\Omega)^{[\gamma,U]}$.

Proof. Let $I:C_0^\infty(\Omega)\rightarrow C^\infty(\Omega)$,
$I(\xi)=\xi$ for $\xi\in C_0^\infty(\Omega)$. Then $I$ is a
continuous linear operator. Hence $U=I^{-1}(V)\in\mathcal
N(C_0^\infty(\Omega))$.

Let $\xi\in C_0^\infty(\Omega)$, $\eta\in U$ and $|t|\leq1$. Then
$\xi\in C^\infty(\Omega)$ and $\eta=I(\eta)\in V$ so
$f_0(\xi+t\eta)=f(\xi+t\eta)=f(\xi)+sf(\eta)=f_0(\xi)+sf_0(\eta)$
where $|s|\leq|\gamma(t)|$. Thus $f_0\in\mathscr
K_{\gamma,U}(C_0^\infty(\Omega),\mathbb C)$.

If $(\xi_\lambda)_{\lambda\in\Delta}$ is a net in
$C_0^\infty(\Omega)$ such that $\xi_\lambda\rightarrow\xi\in
C_0^\infty(\Omega)$. Then $\xi_\lambda=I(\xi_\lambda)\rightarrow
I(\xi)=\xi$ in $C^\infty(\Omega)$ and so
$f_0(\xi_\lambda)=f(\xi_\lambda)\rightarrow f(\xi)=f_0(\xi)$. This
shows that $f_0:C_0^\infty(\Omega)\rightarrow\mathbb C$ is
continuous, i.e., $f_0\in C_0^\infty(\Omega)^{[\gamma,U]}$. \ \
$\square$
\end{Theorem}

\section{Demi-distributions with Compact Support}
\begin{Definition}
Let $f\in\mathscr L_{\gamma,U}(C_0^\infty(\Omega),\mathbb C)$
(resp., $\mathscr L_{\gamma,V}(C^\infty(\Omega),\mathbb C)$) and
$k\in\mathbb N\cup\{0\}$. If for every compact $K\subset\Omega$
there is a $C>0$ such that
$$\big|f(\xi)\big|\leq C\sum_{|\alpha|\leq k}\sup_{x\in K}\big|\partial^\alpha\xi(x)\big|,
\ \forall\,\xi\in C_0^\infty(K)\ (\mbox{resp., }\xi\in
C^\infty(K)),\leqno(1.1)$$ then we say that $f$ is of order $\leq
k$.
\end{Definition}

Let $M\geq1$ and $\gamma(t)=Mt,\ \forall\,t\in\mathbb C$. Then
$\gamma\in C(0)$ and for every $U\in\mathcal N(C_0^\infty(\Omega))$
the family of demi-distributions $C_0^\infty(\Omega)^{[\gamma,U]}$
is a very large extension of $\mathscr D'(\Omega)$
($=C_0^\infty(\Omega)'$), the family of usual distributions (see [2,
Th. 1.1, Th. 2.1]; [1, Th. 1.5, Cor. 1.3]). By Th. 3.6, if $f\in
C_0^\infty(\Omega)^{[\gamma,U]}$ has compact support, then $f$ has
an extension $\widetilde{f}\in C^\infty(\Omega)^{[\gamma,V]}$ and,
conversely, every $f\in C^\infty(\Omega)^{[\gamma,V]}$ has the
restriction $f|_{C_0^\infty(\Omega)}\in
C_0^\infty(\Omega)^{[\gamma,U]}$, where the relations between $U$
and $V$ are very simple.

For $C^\infty(\Omega)^{[\gamma,U]}$ we have a very nice result as
follows.

\begin{Theorem}
Let $M\geq1$, $\gamma(t)=Mt,\ \forall\,t\in\mathbb C,\ V\in\mathcal
N(C^\infty(\Omega))$. Then for every $f\in
C^\infty(\Omega)^{[\gamma,V]}$ there exist compact $L\subset\Omega$,
$C>0$ and $k\in\mathbb N\cup\{0\}$ such that
$$\big|f(\xi)\big|\leq C\sum_{|\alpha|\leq k}\sup_{L}\big|\partial^\alpha\xi\big|,
\ \ \forall\,\xi\in C^\infty(\Omega).\leqno(4.1)$$ Thus, $supp\,f$
is compact, the condition (1.1) holds for $f$, and $f$ is of order
$\leq k$.

Proof. Let $P=\big\{\|\cdot\|_{K,k}:K\mbox{ is a compact subset of
}\Omega,\ k\in\mathbb N\cup\{0\}\big\}$. The topology of
$C^\infty(\Omega)$ is just given by the seminorm family $P$.

There exist $\|\cdot\|_1,\cdots,\|\cdot\|_p\in P$ and
$\varepsilon_1,\cdots,\varepsilon_p\in(0,+\infty)$ such that
$$\bigcap_{j=1}^p\big\{\xi\in C^\infty(\Omega):\|\xi\|_j\leq\varepsilon_j\big\}\subset V.$$
Since $f$ is continuous and $f(0)=0$, there exist
$\|\cdot\|_{p+1},\cdots,\|\cdot\|_m\in P$ and
$\varepsilon_{p+1},\cdots,\varepsilon_m\in(0,+\infty)$ such that
$$\big|f(\xi)\big|<1,\ \
\forall\,\xi\in\bigcap_{j=p+1}^m\big\{\eta\in
C^\infty(\Omega):\|\eta\|_j\leq\varepsilon_j\big\}.$$

Say that $\|\cdot\|_j=\|\cdot\|_{K_j,k_j},\ j=1,2,\cdots,m$, and
$\theta=\min_{1\leq j\leq m}\varepsilon_j$. Then $\theta>0$. Letting
$L=\bigcup_{j=1}^mK_j$, $k=\sum_{j=1}^mk_j$ and, simply,
$\|\cdot\|=\|\cdot\|_{L,k}$, $L$ is compact and $\|\cdot\|\in P$.

If $\xi\in C^\infty(\Omega)$ such that $\|\xi\|\leq\theta$, then
$$\|\xi\|_j=\sum_{|\alpha|\leq k_j}\sup_{K_j}|\partial^\alpha\xi|\leq
\sum_{|\alpha|\leq
k}\sup_L|\partial^\alpha\xi|=\|\xi\|\leq\theta\leq\varepsilon_j,\
j=1,2,\cdots,m.$$ Thus $W=\big\{\xi\in
C^\infty(\Omega):\|\xi\|\leq\theta\big\}\subset\bigcap_{j=1}^m\big\{\xi\in
C^\infty(\Omega):\|\xi\|_j\leq\varepsilon_j\big\}$ and so
$$W\subset V,\ \ \big|f(\xi)\big|<1,\ \ \forall\,\xi\in W.$$

If $\xi\in C^\infty(\Omega)$ such that $\|\xi\|=0$, then
$\|p\xi\|=p\|\xi\|=0$ for all $p\in\mathbb N$ so $p\xi\in W\subset
V$ for all $p\in\mathbb N$ and
$|f(\xi)|=|f(\frac{1}{p}p\xi)|=|s_pf(p\xi)|\leq|s_p|\leq|\gamma(\frac{1}{p})|=M\frac{1}{p}\rightarrow0$
as $p\rightarrow+\infty$. Thus,
$|f(\xi)|=0\leq\frac{M}{\theta}\|\xi\|$.

Let $\xi\in C^\infty(\Omega)$ with $\|\xi\|>0$. Then
$0<\frac{\|\xi\|}{p\theta}\leq1$ for some $p\in\mathbb N$, and
$\|\frac{\theta}{\|\xi\|}\xi\|=\theta$,
$\frac{\theta}{\|\xi\|}\xi\in W\subset V$ so
$|f(\frac{\theta}{\|\xi\|}\xi)|<1$. Hence
\begin{align*}
\big|f(\xi)\big|&=\Big|f\big(p\frac{\|\xi\|}{p\theta}\frac{\theta}{\|\xi\|}\xi\big)\Big|
=\Big|f\big[(p-1)\frac{\|\xi\|}{p\theta}\big(\frac{\theta}{\|\xi\|}\xi\big)+
\frac{\|\xi\|}{p\theta}\big(\frac{\theta}{\|\xi\|}\xi\big)\big]\Big|\\
&=\Big|f\big[(p-1)\frac{\|\xi\|}{p\theta}\big(\frac{\theta}{\|\xi\|}\xi\big)\big]+
s_1f\big(\frac{\theta}{\|\xi\|}\xi\big)\Big|\\
&\qquad\qquad\cdots\\
&=\Big|f\big[\frac{\|\xi\|}{p\theta}\big(\frac{\theta}{\|\xi\|}\xi\big)\big]
+s_{p-1}f\big(\frac{\theta}{\|\xi\|}\xi\big)+\cdots+s_1f\big(\frac{\theta}{\|\xi\|}\xi\big)\Big|\\
&=\Big|\sum_{j=1}^ps_jf\big(\frac{\theta}{\|\xi\|}\xi\big)\Big|\\
&=\Big|\sum_{j=1}^ps_j\Big|\Big|f\big(\frac{\theta}{\|\xi\|}\xi\big)\Big|\\
&\leq\sum_{j=1}^p\big|s_j\big|\leq\sum_{j=1}^p\big|\gamma\big(\frac{\|\xi\|}{p\theta}\big)\big|
=p\big|\gamma(\frac{\|\xi\|}{p\theta})\big|=pM\frac{\|\xi\|}{p\theta}=\frac{M}{\theta}\big\|\xi\big\|.
\end{align*}
Thus we have that
$$\big|f(\xi)\big|\leq\frac{M}{\theta}\|\xi\|=\frac{M}{\theta}\sum_{|\alpha|\leq
k}\sup_L\big|\partial^\alpha\xi\big|,\ \ \forall\,\xi\in
C^\infty(\Omega).$$

If $x\in\Omega\backslash L$, then there is an $\varepsilon>0$ such
that
$N_x=\big\{y\in\Omega:|y-x|=\sqrt{(y_1-x_1)^2+\cdots+(y_n-x_n)^2}\leq\varepsilon\big\}\subset\Omega\backslash
L$ and $|f(\xi)|\leq\frac{M}{\theta}\sum_{|\alpha|\leq
k}\sup_L|\partial^\alpha\xi|=0$ for all $\xi\in C^\infty(N_x)$. Thus
$supp\,f\subset L$ and so $supp\,f$ is compact.

If $K$ is a compact subset of $\Omega$, then for every $\xi\in
C^\infty(K)$ we have that
$$\big|f(\xi)\big|\leq\frac{M}{\theta}\sum_{|\alpha|\leq k}\sup_L\big|\partial^\alpha\xi\big|=
\frac{M}{\theta}\sum_{|\alpha|\leq k}\sup_{L\cap
K}\big|\partial^\alpha\xi\big|
\leq\frac{M}{\theta}\sum_{|\alpha|\leq
k}\sup_K\big|\partial^\alpha\xi\big|,$$ i.e., the condition (1.1)
holds for $f$ and $f$ is of order $\leq k$.\ \ $\square$
\end{Theorem}

Now we can obtain many important facts by the help of Th. 4.1.

\begin{Theorem}
$M\geq1$, $\gamma(t)=Mt$ for $t\in\mathbb C$, $U\in\mathcal
N(C_0^\infty(\Omega))$. If $f\in C_0^\infty(\Omega)^{[\gamma,U]}$
has compact support, then there exist compact $L\subset\Omega$,
$C>0$ and $k\in\mathbb N\cup\{0\}$ such that
$$\big|f(\xi)\big|\leq C\sum_{|\alpha|\leq
k}\sup_L\big|\partial^\alpha\xi\big|,\ \ \forall\,\xi\in
C_0^\infty(\Omega).\leqno(4.1)'$$ Thus, the condition (1.1) holds
for $f$, and $f$ is of order $\leq k$.

Proof. By Th. 3.6 there is a $V\in\mathcal N(C^\infty(\Omega))$ such
that the canonical extension $\widetilde{f}\in
C^\infty(\Omega)^{[\gamma,V]}$. By Th. 4.1, there exist compact
$L\subset\Omega$, $C>0$ and $k\in\mathbb N\cup\{0\}$ such that
$$\big|f(\xi)\big|=\big|\widetilde{f}(\xi)\big|\leq C\sum_{|\alpha|\leq
k}\sup_L\big|\partial^\alpha\xi\big|,\ \ \forall\,\xi\in
C_0^\infty(\Omega).$$

As in the proof of Th. 4.1, (1.1) holds for $f$, and $f$ is of order
$\leq k$.\ \ $\square$
\end{Theorem}

For $\gamma(t)=et\in C(0)$ and $U=\big\{\xi\in C_0^\infty(\mathbb
R):\sup_{0\leq x\leq a}|\xi(x)|\leq1\big\}$ where $a>0$, there
exists demi-distribution $f\in C_0^\infty(\mathbb R)^{(\gamma,U)}$
such that $supp\,f=\{x_0\}$ is compact but the condition (1.1) fails
to hold for $f$ (see Exam. 2.2). However, Th.4.2 shows that if $f\in
C_0^\infty(\Omega)^{[\gamma,U]}$ has compact support then not only
(1.1) holds for $f$ but the more strong (4.1)$'$ holds for $f$.
Thus, the most important properties of demi-distributions heavily
depend on the splitting degree of demi-distributions.

\begin{Theorem}
$M\geq1$, $\gamma(t)=Mt$ for $t\in\mathbb C$. If $f\in
C^\infty(\Omega)^{[\gamma,V]}$ and $f_0(\xi)=f(\xi)$ for $\xi\in
C_0^\infty(\Omega)$, then $f_0\in C_0^\infty(\Omega)^{[\gamma,U]}$
where $U=V\cap C_0^\infty(\Omega)\in\mathcal N(C_0^\infty(\Omega))$,
and there exist compact $L\subset\Omega$, $C>0$ and $k\in\mathbb
N\cup\{0\}$ such that
$$\big|f_0(\xi)\big|\leq C\sum_{|\alpha|\leq k}\sup_L\big|\partial^\alpha\xi\big|,\
\ \forall\,\xi\in C_0^\infty(\Omega)$$ so $supp\,f_0$ is compact,
$supp\,f_0=supp\,f$ and $f_0$ is of order $\leq k$. Moreover,
$f=\widetilde{f_0}$, the canonical extension of $f_0$.

Proof. By Th. 3.7 and Th. 4.1, we only need to show
$\widetilde{f_0}=f$.

By Th. 4.1, there exist compact $L\subset\Omega$, $C>0$ and
$k\in\mathbb N\cup\{0\}$ such that
$$\big|f(\xi)\big|\leq C\sum_{|\alpha|\leq k}\sup_L\big|\partial^\alpha\xi\big|,\
\ \forall\,\xi\in C^\infty(\Omega).$$ If $x\in\Omega\backslash
supp\,f$ then there is an open $N_x\subset\Omega\backslash supp\,f$
such that $f(\xi)=0$ for all $\xi\in C^\infty(N_x)$ and so
$f_0(\xi)=f(\xi)=0,\ \forall\,\xi\in C_0^\infty(N_x)\subset
C^\infty(N_x)$, that is, $x\not\in supp\,f_0$. Thus
$supp\,f_0\subset supp\, f\subset L$.

Pick a $\mathcal X\in C_0^\infty(\Omega)$ such that $\mathcal X=1$
in a neighborhood of $L$. Then by Th. 3.4 and Def. 3.2 we have that
$\widetilde{f_0}(\xi)=f_0(\mathcal X\xi)=f(\mathcal X\xi),\
\forall\,\xi\in C^\infty(\Omega)$.

Let $\xi\in C^\infty(\Omega)$ and pick a $p\in\mathbb N$ such that
$\frac{1}{p}(1-\mathcal X)\xi\in V$. Since $1-\mathcal X=0$ in a
neighborhood of $L$, $\partial^\alpha[\frac{1}{p}(1-\mathcal
X)\xi](x)=0$ for all $x\in L$ and all multi-index $\alpha$. Then
$$\big|f\big[\frac{1}{p}(1-\mathcal X)\xi\big]\big|\leq
C\sum_{|\alpha|\leq
k}\sup_L\big|\partial^\alpha\big[\frac{1}{p}(1-\mathcal
X)\xi\big]\big|=0 \mbox{, i.e., }f\big[\frac{1}{p}(1-\mathcal
X)\xi\big]=0,$$
\begin{align*}
f(\xi)=&f\big[\mathcal X\xi+(1-\mathcal X)\xi\big]=f\big[\mathcal
X\xi+(p-1)\frac{1}{p}(1-\mathcal X)\xi+\frac{1}{p}(1-\mathcal
X)\xi\big]\\
&=f\big[\mathcal X\xi+(p-1)\frac{1}{p}(1-\mathcal
X)\xi\big]+sf\big[\frac{1}{p}(1-\mathcal X)\xi\big]\\
&=f\big[\mathcal X\xi+(p-1)\frac{1}{p}(1-\mathcal X)\xi\big]\\
&\qquad\qquad\cdots\\
&=f(\mathcal X\xi)=f_0(\mathcal X\xi)=\widetilde{f_0}(\xi).
\end{align*}

Thus $f=\widetilde{f_0}$ and
$supp\,f_0=supp\,\widetilde{f_0}=supp\,f$ by Th. 3.6.\ \ $\square$
\end{Theorem}

\begin{Corollary}
$M\geq1$, $\gamma(t)=Mt,\ \forall\,t\in\mathbb C$. Then
$$\bigcup_{V\in\mathcal N(C^\infty(\Omega))}C^\infty(\Omega)^{[\gamma,V]}=
\bigcup_{U\in\mathcal N(C_0^\infty(\Omega))}\big\{\widetilde{f}:f\in
C_0^\infty(\Omega)^{[\gamma,U]},\ supp\,f\mbox{ is compact}\big\}.$$
\end{Corollary}

Now we can improve Th. 2.4 as follows.

\begin{Corollary}
$M\geq1$, $\gamma(t)=Mt$ for $t\in\mathbb C$. Let $f\in\mathscr
K_{\gamma,U}(C_0^\infty(\Omega),\mathbb C)$ for which $supp\,f$ is
compact. Then $f$ is continuous if and only if the condition (1.1)
holds for $f$.

Proof. If $f$ is continuous, then $f\in
C_0^\infty(\Omega)^{[\gamma,U]}$ and Th. 4.2 shows that (1.1) holds
for $f$.

Conversely, suppose that (1.1) holds for $f$. Since $supp\,f$ is
compact, there is a $\mathcal X\in C_0^\infty(\Omega)$ such that
$\mathcal X=1$ in a neighborhood of $supp\,f$. Then
$L=supp\,\mathcal X$ is compact. By Lemma 3.4, $V=\big\{\xi\in
C^\infty(\Omega):\mathcal X\xi\in U\big\}\in\mathcal
N(C^\infty(\Omega))$ and the canonical extension
$\widetilde{f}\in\mathscr K_{\gamma,V}(C^\infty(\Omega),\mathbb C)$.
In fact, for $\xi\in C^\infty(\Omega)$, $\eta\in V$ and $|t|\leq1$,
$\widetilde{f}(\xi+t\eta)=f(\mathcal X\xi+t\mathcal
X\eta)=f(\mathcal X\xi)+sf(\mathcal
X\eta)=\widetilde{f}(\xi)+s\widetilde{f}(\eta)$ where
$|s|\leq|\gamma(t)|$.

Since (1.1) holds for $f$, $f$ is sequentially continuous by Th.
2.3. If $\xi_v\rightarrow0$ in $C^\infty(\Omega)$ then $\mathcal
X\xi_v\rightarrow0$ in $C_0^\infty(\Omega)$ by Lemma 3.4, and
$\widetilde{f}(\xi_v)=f(\mathcal X\xi_v)\rightarrow
f(0)=0=\widetilde{f}(0)$, that is, $\widetilde{f}$ is continuous at
$0\in C^\infty(\Omega)$ since $C^\infty(\Omega)$ is a Fr\'echet
space. Then $\widetilde{f}$ is continuous by Lemma 3.5,
$\widetilde{f}\in C^\infty(\Omega)^{[\gamma,V]}$. By Th. 3.7, $f$ is
continuous.\ \ $\square$
\end{Corollary}

For $x=(x_1,\cdots,x_n)\in\mathbb R^n$,\ \
$|x|=\sqrt{x_1^2+\cdots+x_n^2}$.

\begin{Lemma}
Let $K$ and $F$ be nonempty subsets of $\Omega$. If $K$ is compact
and $F$ is closed in $\Omega$ and $K\cap F=\emptyset$, then there
exist $x_0\in K$ and $y_0\in\overline{F}$ such that $\inf_{x\in
K,y\in F}|x-y|=|x_0-y_0|>0$.

Proof. Let $d=\inf_{x\in K,y\in F}|x-y|$. There exist sequences
$\{x_v\}\subset K$ and $\{y_v\}\subset F$ such that
$d=\lim_v|x_v-y_v|$. Since $K$ is compact and $\{x_v-y_v\}$ is
bounded, we may assume that $x_v\rightarrow x_0\in K$ and
$x_v-y_v\rightarrow b\in\mathbb R^n$. Then
$y_v=y_v-x_v+x_v\rightarrow x_0-b=y_0$ and
$|x_v-y_v|\rightarrow|x_0-y_0|$, $d=|x_0-y_0|$.

If $y_0\in F$ then $y_0\not\in K$ so $y_0\neq x_0$ and
$d=|x_0-y_0|>0$. If $y_0\not\in F$ then $y_0\not\in\Omega$ so
$y_0\neq x_0$ and $d=|x_0-y_0|>0$.\ \ $\square$
\end{Lemma}

We have a fact which is different from Lemma 2.2 as follows.

\begin{Theorem}
$M\geq1$, $\gamma(t)=Mt$ for $t\in\mathbb C$. If $f\in
C^\infty(\Omega)^{[\gamma,V]}$ and $\xi\in C^\infty(\Omega)$ such
that $supp\,f\cap supp\,\xi=\emptyset$, then $f(\xi)=0$.

Proof. By Th. 4.1, $supp\,f$ is compact. Then $\inf\big\{|x-y|:x\in
supp\,f,\ y\in supp\,\xi\big\}=d>0$ by Lemma 4.1.

Let $f_0(\xi)=f(\xi)$ for $\xi\in C_0^\infty(\Omega)$ and $U=V\cap
C_0^\infty(\Omega)$. By Th. 4.3, $f_0\in
C_0^\infty(\Omega)^{[\gamma,U]}$ and the canonical extension
$\widetilde{f_0}=f$, $supp\,f_0=supp\,f$. Since $d>0$, for
$\varepsilon\in(0,d/3)$ there is a $\mathcal X\in
C_0^\infty(\Omega)$ such that $0\leq\mathcal X\leq1$ and
$$\mathcal X=1\mbox{ in }G=\big\{y\in\Omega:|y-x|<\varepsilon\mbox{ for some }x\in supp\,f_0\big\},$$
$$\mathcal X=0\mbox{ outside }B_{3\varepsilon}=
\big\{y\in\Omega:|y-x|\leq3\varepsilon\mbox{ for some }x\in supp\,
f_0\big\}.$$ Then $supp\,\mathcal X\subset B_{3\varepsilon}$,
$supp\,f_0\cap supp\,(\mathcal X\xi)\subset supp\,f\cap
supp\,\xi=\emptyset$ and $supp\,[(1-\mathcal X)\xi]\cap
supp\,f_0\subset(\Omega\backslash G)\cap supp\,f_0=\emptyset$. Hence
$\xi=\mathcal X\xi+(1-\mathcal X)\xi$ is a $f_0$-decomposition of
$\xi$ and $supp\,(\mathcal X\xi)\cap supp\,f_0=\emptyset$. Then
$f(\xi)=\widetilde{f_0}(\xi)=f_0(\mathcal X\xi)=0$ by Lemma 2.2.\ \
$\square$
\end{Theorem}

\begin{Theorem}
$M\geq1$, $\gamma(t)=Mt$ for $t\in\mathbb C$. If $f\in
C^\infty(\Omega)^{[\gamma,V]}$ is of order $\leq k$ and $\xi\in
C^\infty(\Omega)$ such that $\partial^\alpha\xi(x)=0$ when
$|\alpha|\leq k$ and $x\in supp\,f$, then $f(\xi)=0$.

Proof. By Th. 4.1, $supp\,f$ is compact so for sufficiently small
$\varepsilon>0$ the set $B_\varepsilon=\big\{y\in\mathbb
R^n:|y-x|\leq\varepsilon$ for some $x\in supp\,f\big\}$ is compact
and contained in $\Omega$. There is a $\mathcal X_\varepsilon\in
C_0^\infty(\Omega)$ with $0\leq\mathcal X_\varepsilon\leq1$ such
that $\mathcal X_\varepsilon=1$ in a neighborhood of $supp\,f$ and
$\mathcal X_\varepsilon=0$ outside $B_\varepsilon$ [4, p.46].
Moreover, $|\partial^\alpha\mathcal X_\varepsilon|\leq C_\alpha
\varepsilon^{-|\alpha|}$ where $C_\alpha$ is independent of
$\varepsilon$ [4, p.5] so there is a $C>0$ such that
$|\partial^\alpha\mathcal X_\varepsilon|\leq
C\varepsilon^{-|\alpha|}$ for all $|\alpha|\leq k$ and all
$\varepsilon\in(0,\varepsilon_0]$, where $\varepsilon_0>0$ and
$B_{\varepsilon_0}\subset\Omega$.

For $\varepsilon\in(0,\varepsilon_0]$ pick a $p\in\mathbb N$ for
which $\frac{1}{p}(1-\mathcal X_\varepsilon)\xi\in V$. Since
$supp\,[\frac{1}{p}(1-\mathcal X_\varepsilon)\xi]\cap
supp\,f=\emptyset$, $f[\frac{1}{p}(1-\mathcal X_\varepsilon)\xi]=0$
by Th. 4.4. Then
$$f(\xi)=f\big[\xi\mathcal X_\varepsilon+p\frac{1}{p}(1-\mathcal X_\varepsilon)\xi\big]=
f\big[\xi\mathcal X_\varepsilon+(p-1)\frac{1}{p}(1-\mathcal
X_\varepsilon)\xi\big]=\cdots=f(\xi\mathcal X_\varepsilon).$$

Since $f$ is of order $\leq k$, there is an $A>0$ such that
$$\big|f(\eta)\big|\leq A\sum_{|\alpha|\leq
k}\sup_{B_\varepsilon}\big|\partial^\alpha\eta\big|,\ \
\forall\,\eta\in C^\infty(B_\varepsilon).$$ If
$0<\varepsilon'<\varepsilon$ and $\eta\in
C^\infty(B_{\varepsilon'})$, then $\eta\in C^\infty(B_\varepsilon)$
and
$\sup_{B_\varepsilon'}|\partial^\alpha\eta|=\sup_{B_\varepsilon}|\partial^\alpha\eta|$
for all $\alpha$ so $|f(\eta)|\leq A\sum_{|\alpha|\leq
k}\sup_{B_\varepsilon}|\partial^\alpha\eta|= A\sum_{|\alpha|\leq
k}\sup_{B_\varepsilon'}|\partial^\alpha\eta|$. Hence the constant
$A$ is available for all $\varepsilon'\in(0,\varepsilon]$. Since
$\xi\mathcal X_\varepsilon\in C_0^\infty(B_\varepsilon)=
C^\infty(B_\varepsilon)$, it follows from the Leibniz formula that
\begin{align*}
\big|f(\xi)\big|=\big|f(\xi\mathcal X_\varepsilon)\big|&\leq
A\sum_{|\alpha|\leq
k}\sup_{B_\varepsilon}\big|\partial^\alpha(\xi\mathcal
X_\varepsilon)\big|\leq A_1\sum_{|\alpha|+|\beta|\leq
k}\sup_{B_\varepsilon}\big|\partial^\alpha\xi\big|\big|\partial^\beta\mathcal
X_\varepsilon\big|\\
&\leq A_1 C\sum_{|\alpha|\leq
k}\varepsilon^{|\alpha|-k}\sup_{B_\varepsilon}\big|\partial^\alpha\xi\big|,
\end{align*}
where both $A_1$ and $C$ are independent of $\varepsilon$. Observing
$\partial^\alpha\xi(x)=0$ for all $x\in supp\,f$ and $|\alpha|\leq
k$, we have
$\lim_{\varepsilon\rightarrow0}\varepsilon^{|\alpha|-k}\sup_{B_\varepsilon}|\partial^\alpha\xi|=0$
for all $|\alpha|\leq k$ [4, p.46]. Thus $f(\xi)=0$.\ \ $\square$
\end{Theorem}

Notice that in the notation $(\bullet-a)^\alpha$ the symbol
$\bullet$ denotes the variable, that is, $(\bullet-a)^\alpha$ is a
function such that
$[(\bullet-a)^\alpha](x)=(x-a)^\alpha=(x_1-a_1)^{\alpha_1}\cdots(x_n-a_n)^\alpha_n$
[4, p.47].

\begin{Corollary}
$M\geq1$, $\gamma(t)=Mt$ for $t\in\mathbb C$. If $f\in
C^\infty(\Omega)^{[\gamma,V]}$ is of order $k$ and $supp\,f=\{y\}$,
a singleton, then we have
$$f(\xi)=f\big[\sum_{|\alpha|\leq k}\partial ^\alpha\xi(y)(\bullet-y)^\alpha/(\alpha!)\big],\ \
\forall\,\xi\in C^\infty(\Omega),$$ where for $0=(0,\cdots,0)$ the
term $\alpha^0\xi(y)(\bullet-y)^0/(0!) =\xi(y)$ is the function
$\eta\in C^\infty(\Omega)$ for which $\eta(x)=\xi(y),\
\forall\,x\in\Omega$.

Proof. Let $\xi\in C^\infty(\Omega)$. We expand $\xi$ in a Taylor
series
$$\xi(x)=\sum_{|\alpha|\leq k}\partial^\alpha\xi(y)(x-y)^\alpha/(\alpha!)+\psi(x).$$
Pick an $m\in\mathbb N$ for which $\frac{1}{m}\psi\in V$. Since
$\partial^\alpha(\frac{1}{m}\psi)(y)=\frac{1}{m}\partial^\alpha\psi(y)=0$
when $|\alpha|\leq k$, $f(\frac{1}{m}\psi)=0$ by Th. 4.5. Hence
\begin{align*}
f(\xi)&=f\big[\sum_{|\alpha|\leq
k}\partial^\alpha\xi(y)(\bullet-y)^\alpha/(\alpha!)+(m-1)\frac{1}{m}\psi+\frac{1}{m}\psi\big]\\
&=f\big[\sum_{|\alpha|\leq
k}\partial^\alpha\xi(y)(\bullet-y)^\alpha/(\alpha!)+(m-1)\frac{1}{m}\psi\big]\\
&\qquad\qquad\cdots\\
&=f\big[\sum_{|\alpha|\leq
k}\partial^\alpha\xi(y)(\bullet-y)^\alpha/(\alpha!)\big].\ \ \square
\end{align*}
\end{Corollary}

Even the simplest case is somewhat of an interesting observation.

\begin{Example}
Let $y\in\Omega$ and $f(\xi)=\sin|\xi(y)|$ for $\xi\in
C^\infty(\Omega)$. Clearly, $V=\big\{\xi\in
C^\infty(\Omega):|\xi(y)|\leq1\big\}\in\mathcal
N(C^\infty(\Omega))$. If $\xi\in C^\infty(\Omega)$, $\eta\in V$ and
$|t|\leq1$, then
$f(\xi+t\eta)=\sin|\xi(y)+t\eta(y)|=\sin(|\xi(y)|+s|\eta(y)|)=\sin|\xi(y)|+\theta\sin|\eta(y)|=f(\xi)+\theta
f(\eta)$ where $|\theta|\leq\frac{\pi}{2}|s|\leq\frac{\pi}{2}|t|$.
Letting $\gamma(t)=\frac{\pi}{2}t$ for $t\in\mathbb C$, $f\in
C^\infty(\Omega)^{[\gamma,V]}$ and $supp\,f=\{y\}$. For every
compact $K\subset\Omega$ and $\xi\in C^\infty(K)=C_0^\infty(K)$ we
have that
$$\big|f(\xi)\big|=\big|\sin|\xi(y)|\big|=0\leq\sup_K\big|\partial^0\xi\big|\mbox{ when }y\not\in K,$$
$$\big|f(\xi)\big|=\big|\sin|\xi(y)|\big|\leq\big|\xi(y)\big|\leq\sup_K|\xi|=\sup_K\big|\partial^0\xi\big|
\mbox{ when }y\in K.$$ Thus, $f$ is of order $0$.
\end{Example}

\begin{Corollary}
$M\geq1$, $\gamma(t)=Mt$ for $t\in\mathbb C$. If $f\in
C^\infty(\Omega)^{[\gamma,V]}$ is of order $0$ and $supp\,f=\{y\}$,
then
$$f(\xi)=f\big(\xi(y)\big),\ \ \forall\,\xi\in C^\infty(\Omega),$$
where $\xi(y)$ is a function in $C^\infty(\Omega)$ such that
$\xi(y)(x)=\xi(y),\ \forall\,x\in\Omega$.
\end{Corollary}

In [2] we gave a very clear-cut characterization of demi-linear
functions in $\mathscr K_{M,\varepsilon}(\mathbb R,\mathbb R)$ [2,
Th. 1.1]. We have a similar clear-cut description for demi-linear
functions in $\mathscr K_{M,\varepsilon}(\mathbb C,\mathbb C)$ as
follows.

\begin{Lemma}
Let $g:\mathbb C\rightarrow\mathbb C$ be a function such that
$g(0)=0$ and $g'(z_0)\neq0$ for some $z_0\in\mathbb C$. Let
$\varepsilon>0$. Then $g\in\mathscr K_{M,\varepsilon}(\mathbb
C,\mathbb C)$ for some $M\geq1$ if and only if

(1) $g$ is continuous,

(2) $g(z)\neq0$ for $0<|z|\leq\varepsilon$,

(3) $\inf_{0<|u|\leq\varepsilon}\big|\frac{g(u)}{u}\big|>0$,

(4) $\sup_{z,u\in\mathbb
C,\,0<|u|\leq\varepsilon}\big|\frac{g(z+u)-g(z)}{u}\big|<+\infty$.

Proof. Suppose that $g\in \mathscr K_{M,\varepsilon}(\mathbb
C,\mathbb C)$ where $M\geq 1$. If $z_k\rightarrow z$ in $\mathbb C$,
then for sufficiently large $k\in\mathbb N$ we have that
$|\frac{z_k-z}{\varepsilon}|<1$ and
$g(z_k)=g(z+\frac{z_k-z}{\varepsilon}\varepsilon)=g(z)+s_kg(\varepsilon)$
where $|s_k|\leq M|\frac{z_k-z}{\varepsilon}|\rightarrow0$ so
$g(z_k)\rightarrow g(z)$, that is, $g$ is continuous. Assume that
$g(u)=0$ for some $0<|u|\leq\varepsilon$ and $z\in\mathbb C$,
$z\neq0$. Then $|\frac{z}{ku}|<1$ for some $k\in\mathbb N$ and
$g(z)=g(k\frac{z}{ku}u)=g[(k-1)\frac{z}{ku}u+\frac{z}{ku}u]=
g[(k-1)\frac{z}{ku}u]+s_1g(u)=g[(k-1)\frac{z}{ku}u]=\cdots=g(\frac{z}{ku}u)=s_kg(u)=0$.
Thus $g=0$ but $g'(z_0)\neq0$. This contradiction shows that (2)
holds for $g$.

If $\,\inf_{0<|u|\leq\varepsilon}\big|\frac{g(u)}{u}\big|=0$, then
$\frac{g(u_k)}{u_k}\rightarrow 0$ for some
$\{u_k\}\subset\big\{z\in\mathbb
C:|z|\leq\varepsilon\big\}\setminus\{0\}$. May assume that
$u_k\rightarrow u_0$. If $u_0\neq0$, then
$|\frac{g(u_k)}{u_k}|\rightarrow |\frac{g(u_0)}{u_0}|>0$ by (1) and
(2), a contradiction. So $u_0=0$, $u_k\rightarrow0$. Then
$g(z_0+u_k)=g(z_0)+s_kg(u_k)$ where $|s_k|\leq M|1|=M$, and $0\neq
g'(z_0)=\lim_k\frac{g(z_0+u_k)-g(z_0)}{u_k}=\lim_k\frac{s_kg(u_k)}{u_k}=0$.
This contradiction shows that (3) holds for $g$.

Let $z,u\in \mathbb C$, $0<|u|\leq\varepsilon$. Since
$g(u)=g\big(\frac{u}{\varepsilon}\varepsilon\big)=sg(\varepsilon)$
where $|s|\leq
M\big|\frac{u}{\varepsilon}\big|=\frac{M}{\varepsilon}|u|$ and
$g(z+u)=g(z)+s_1g(u)$ where $|s_1|\leq M|1|=M$,
$|g(u)|\leq\frac{M}{\varepsilon}|g(\varepsilon)u|$ and
$\Big|\frac{g(z+u)-g(z)}{u}\Big|=\Big|\frac{s_1g(u)}{u}\Big|
\leq\frac{M^2}{\varepsilon}|g(\varepsilon)|$. Thus, (4) holds for
$g$.

Conversely, assume that (1), (2), (3) and (4) hold for $g$. Since
$g(0)=0$ and
$\inf_{0<|u|\leq\varepsilon}\big|\frac{g(u)}{u}\big|=\inf_{0<|u|\leq\varepsilon}\big|\frac{g(0+u)-g(0)}{u}\big|\leq
\sup_{0<|u|\leq\varepsilon}\big|\frac{g(0+u)-g(0)}{u}\big|$,
$$M=\Big[\sup_{z,u\in \mathbb
C,\,0<|u|\leq\varepsilon}\big|\frac{g(z+u)-g(z)}{u}\big|\Big]/\inf_{0<|u|\leq\varepsilon}\big|\frac{g(u)}{u}\big|\geq1.$$
Let $z,u,t\in\mathbb C$, $0<|u|\leq\varepsilon$, $0<|t|\leq1$. Then
$g(u)\neq0$ by (2) and
$g(z+tu)\linebreak=g(z)+g(z+tu)-g(z)=g(z)+\Big[\frac{g(z+tu)-g(z)}{tu}\frac{u}{g(u)}t\Big]g(u)$,
where
$$\Big|\frac{g(z+tu)-g(z)}{tu}\frac{u}{g(u)}t\Big|=\Big|\frac{g(z+tu)-g(z)}{tu}/\frac{g(u)}{u}\Big|\big|t\big|\leq
M\big|t\big|.$$ Thus, $g\in\mathscr K_{M,\varepsilon}(\mathbb
C,\mathbb C)$.\ \ $\square$
\end{Lemma}

We now have the following representation theorem.

\begin{Theorem}
$M\geq1$, $\gamma(t)=Mt$ for $t\in\mathbb C$. If $f\in
C^\infty(\Omega)^{[\gamma,V]}$ is of order $0$ and $supp\,f=\{y\}$,
then there exist $\varepsilon>0$ and $g\in\mathscr
K_{M,\varepsilon}(\mathbb C,\mathbb C)$ such that
$$f(\xi)=g\big(\xi(y)\big),\ \ \forall\,\xi\in C^\infty(\Omega).\leqno(4.2)$$
Conversely, every $\varepsilon>0$ and $g\in\mathscr
K_{M,\varepsilon}(\mathbb C,\mathbb C)$ give a $V\in\mathcal
N(C^\infty(\Omega))$ and a $f\in C^\infty(\Omega)^{[\gamma,V]}$
through (4.2) such that $f$ is of order $0$ and $supp\,f=\{y\}$.

Proof. Suppose that $f\in C^\infty(\Omega)^{[\gamma,V]}$ is of order
$0$ and $supp\,f=\{y\}$. For $z\in\mathbb C$ let $\zeta_z(x)=z$ for
all $x\in\Omega$. Then $\zeta_z\in C^\infty(\Omega)$ and
$\lim_{z\rightarrow0}\zeta_z=0$ in $C^\infty(\Omega)$. Hence there
is an $\varepsilon>0$ such that $\zeta_z\in V$ when
$|z|\leq\varepsilon$.

Define $g:\mathbb C\rightarrow\mathbb C$ by $g(z)=f(\zeta_z),\
\forall\,z\in\mathbb C$. Then $g(0)=f(\zeta_0)=f(0)=0$. For
$z,u,t\in\mathbb C$ with $|u|\leq\varepsilon$ and $|t|\leq1$,
$\zeta_u\in V$ and
$$g(z+tu)=f(\zeta_{z+tu})=f(\zeta_z+\zeta_{tu})=f(\zeta_z+t\zeta_u)=f(\zeta_z)+sf(\zeta_u)=g(z)+sg(u),$$
where $|s|\leq|\gamma(t)|=M|t|$. Thus $g\in\mathscr
K_{M,\varepsilon}(\mathbb C,\mathbb C)$. By Cor. 4.4 we have
$$f(\xi)=f\big(\xi(y)\big)=f\big(\zeta_{\xi(y)}\big)=g\big(\xi(y)\big),\ \ \forall\,\xi\in C^\infty(\Omega).$$

Conversely, let $\varepsilon>0$, $g\in\mathscr
K_{M,\varepsilon}(\mathbb C,\mathbb C)$ and $y\in\Omega$. Since the
Dirac measure $\delta_y:C^\infty(\Omega)\rightarrow\mathbb C$,
$\delta_y(\xi)=\xi(y)$ is continuous and $\delta_y(0)=0$,
$V=\big\{\xi\in
C^\infty(\Omega):|\xi(y)|\leq\varepsilon\big\}=\delta_y^{-1}([-\varepsilon,\varepsilon])\in\mathcal
N(C^\infty(\Omega))$. Then define
$f:C^\infty(\Omega)\rightarrow\mathbb C$ by $f(\xi)=g(\xi(y)),\
\xi\in C^\infty(\Omega)$. For $\xi\in C^\infty(\Omega)$, $\eta\in V$
and $|t|\leq1$, $|\eta(y)|\leq\varepsilon$ and
$f(\xi+t\eta)=g((\xi+t\eta)(y))=g(\xi(y)+t\eta(y))=g(\xi(y))+sg(\eta(y))=f(\xi)+sf(\eta)$
where $|s|\leq|\gamma(t)|=M|t|$. Thus $f\in\mathscr
K_{\gamma,V}(C^\infty(\Omega),\mathbb C)$ and $f$ is continuous
because $\xi_\lambda\rightarrow\xi$ in $C^\infty(\Omega)$ implies
$\xi_\lambda(y)\rightarrow\xi(y)$ and
$f(\xi_\lambda)=g(\xi_\lambda(y))\rightarrow g(\xi(y))=f(\xi)$ by
Lemma 4.2, that is, $f\in C^\infty(\Omega)^{[\gamma,V]}$.

Let $y_0\in\Omega$, $y_0\neq y$. Pick a $\theta>0$ such that
$K=\big\{x\in\mathbb
R^n:|x-y_0|\leq\theta\big\}\subset\Omega\backslash\{y\}$. Then for
every $\xi\in C^\infty(K)$ we have $\xi(y)=0$ and
$f(\xi)=g(\xi(y))=g(0)=0$. This shows that $y_0\not\in supp\,f$,
$supp\,f\subset\{y\}$. If $G$ is an open set in $\mathbb R^n$ such
that $y\in G\subset\Omega$, then there is a $\mathcal X\in
C_0^\infty(G)$ such that $0<|\mathcal X(y)|\leq\varepsilon$. By
Lemma 4.2, $f(\mathcal X)=g(\mathcal X(y))\neq0$ so $y\in supp\,f$
and $supp\,f=\{y\}$.

Let $K$ be a compact subset of $\Omega$. If $y\in K$ and $\xi\in
C^\infty(K)$ then there is a $p\in\mathbb N$ such that
$|\frac{\xi(y)}{p\varepsilon}|<1$ and
\begin{align*}
\big|f(\xi)\big|&=\big|g\big(\xi(y)\big)\big|=\big|g\big(p\frac{\xi(y)}{p\varepsilon}\varepsilon\big)\big|=
\big|g\big[(p-1)\frac{\xi(y)}{p\varepsilon}\varepsilon+\frac{\xi(y)}{p\varepsilon}\varepsilon\big]\big|\\
&=\big|g\big[(p-1)\frac{\xi(y)}{p\varepsilon}\varepsilon\big]+s_1g(\varepsilon)\big|\\
&\qquad\quad\cdots\\
&=\big|g(\frac{\xi(y)}{p\varepsilon}\varepsilon)+s_{p-1}g(\varepsilon)+\cdots+s_1g(\varepsilon)\big|\\
&=\big|s_pg(\varepsilon)+s_{p-1}g(\varepsilon)+\cdots+s_1g(\varepsilon)\big|
=\big|\sum_{v=1}^ps_v\big|\big|g(\varepsilon)\big|,
\end{align*}
where each
$|s_v|\leq|\gamma(\frac{\xi(y)}{p\varepsilon})|=M\frac{|\xi(y)|}{p\varepsilon}$
so $|\sum_{v=1}^ps_v|\leq\sum_{v=1}^p|s_v|\leq
pM\frac{|\xi(y)|}{p\varepsilon}=\frac{M}{\varepsilon}|\xi(y)|$.
Therefore,
$$\big|f(\xi)\big|\leq\frac{M}{\varepsilon}\big|g(\varepsilon)\big|\big|\xi(y)\big|\leq
\frac{M}{\varepsilon}\big|g(\varepsilon)\big|\sup_K\big|\xi\big|
=\frac{M}{\varepsilon}\big|g(\varepsilon)\big|\sup_K\big|\partial^0\xi\big|.$$
If $y\not\in K$ and $\xi\in C^\infty(K)$ then $\xi(y)=0$ and
$$\big|f(\xi)\big|=\big|g(\xi(y))\big|=\big|g(0)\big|=0
\leq\frac{M}{\varepsilon}\big|g(\varepsilon)\big|\sup_K\big|\partial^0\xi\big|.$$

Thus, $f$ is of order $0$. Moreover, the constant
$\frac{M}{\varepsilon}|g(\varepsilon)|$ is available for all compact
$K\subset\Omega$.\ \ $\square$
\end{Theorem}

\begin{Corollary}
$M\geq1$, $\gamma(t)=Mt$ for $t\in\mathbb C$. If $f\in
C^\infty(\Omega)^{[\gamma,V]}$ is of order $0$ and $supp\,f=\{y\}$,
then there exist $\varepsilon>0$ and $g\in\mathscr
K_{M,\varepsilon}(\mathbb C,\mathbb C)$ such that
$$f(\xi)=a_\xi M\frac{g(\varepsilon)}{\varepsilon}\xi(y),\ \forall\,\xi\in C^\infty(\Omega),$$
where $|a_\xi|\leq1$. Hence,
$A=M|\frac{g(\varepsilon)}{\varepsilon}|>0$ and $|f(\xi)|\leq
A|\xi(y)|,\ \forall\,\xi\in C^\infty(\Omega)$.

Proof. By Th. 4.6, there exist $\varepsilon>0$ and $g\in\mathscr
K_{M,\varepsilon}(\mathbb C,\mathbb C)$ such that
$$f(\xi)=g\big(\xi(y)\big),\ \ \forall\,\xi\in C^\infty(\Omega).$$
Let $\xi\in C^\infty(\Omega)$. Pick a $p\in\mathbb N$ such that
$|\frac{\xi(y)}{p\varepsilon}|<1$. Then
$$g\big(\xi(y)\big)=g\big(p\frac{\xi(y)}{p\varepsilon}\varepsilon\big)=
g\big[(p-1)\frac{\xi(y)}{p\varepsilon}\varepsilon\big]+s_1g(\varepsilon)=\cdots=\big(\sum_{v=1}^ps_v\big)g(\varepsilon),$$
where $|\sum_{v=1}^ps_v|\leq\frac{M}{\varepsilon}|\xi(y)|$, i.e.,
$|g(\xi(y))|=|\sum_{v=1}^ps_v||g(\varepsilon)|\leq\frac{M}{\varepsilon}|g(\varepsilon)\xi(y)|$.
Hence there is an $a_\xi\in\mathbb C$ with $|a_\xi|\leq1$ such that
$g(\xi(y))=a_\xi M\frac{g(\varepsilon)}{\varepsilon}\xi(y)$.

Since $supp\,f=\{y\}\neq\emptyset$, $f\neq0$ and so $g\neq0$. By
Lemma 4.2, $A=M|\frac{g(\varepsilon)}{\varepsilon}|>0$.\ \ $\square$
\end{Corollary}

\begin{Corollary}
$M\geq1$, $\gamma(t)=Mt$ for $t\in\mathbb C$. If $f\in
C^\infty(\Omega)^{[\gamma,V]}$ is of order $0$ and $supp\,f=\{y\}$,
then $f$ is Lipschitz, that is, there is a $A>0$ such that
$$\big|f(\xi)-f(\eta)\big|\leq A\big|\xi(y)-\eta(y)\big|,\ \ \forall\,\xi,\eta\in C^\infty(\Omega),$$
$$\big|f(\xi)\big|\leq A\big|\xi(y)\big|,\ \ \forall\,\xi\in C^\infty(\Omega).$$

Proof. By Th. 4.6, there exist $\varepsilon>0$ and $g\in\mathscr
K_{M,\varepsilon}(\mathbb C,\mathbb C)$ such that
$$f(\xi)=g\big(\xi(y)\big),\ \ \forall\,\xi\in C^\infty(\Omega).$$
Since $supp\,f=\{y\}\neq\emptyset$, $f\neq0$ and so $g\neq0$.

Let $z,u\in\mathbb C$, $0<|u|\leq\varepsilon$. Since $g(0)=0$,
$g(u)=g(\frac{u}{\varepsilon}\varepsilon)=g(0+\frac{u}{\varepsilon}\varepsilon)=g(0)+sg(\varepsilon)=sg(\varepsilon)$
where $|s|\leq M|\frac{u}{\varepsilon}|=\frac{M}{\varepsilon}|u|$
and so
$|g(u)|=|sg(\varepsilon)|\leq\frac{M}{\varepsilon}|ug(\varepsilon)|$.
Then $g(z+u)=g(z)+s_1g(u)$ where $|s_1|\leq M|1|=M$ and
$$\big|g(z+u)-g(z)\big|=\big|s_1g(u)\big|\leq\frac{M^2}{\varepsilon}\big|g(\varepsilon)\big|\big|u\big|.$$
Letting $A=\frac{M^2}{\varepsilon}|g(\varepsilon)|$ we have that
$$\big|g(z+u)-g(z)\big|\leq A|u|,\
\ \forall\,z\in\mathbb C,\ u\in\{z\in\mathbb
C:|z|\leq\varepsilon\}.\leqno(4.3)$$

Let $\xi,\eta\in C^\infty(\Omega)$. If $\xi(y)=\eta(y)$, then
$$\big|f(\xi)-f(\eta)\big|=\big|g(\xi(y))-g(\eta(y))\big|=0=A\big|\xi(y)-\eta(y)\big|.$$
If $\xi(y)\neq\eta(y)$, then
$0<\frac{1}{p}|\xi(y)-\eta(y)|<\varepsilon$ for some $p\in\mathbb
N$. Let $u=\frac{1}{p}[\xi(y)-\eta(y)]$, then $\xi(y)=\eta(y)+pu$
and
$$f(\xi)-f(\eta)=g\big(\xi(y)\big)-g\big(\eta(y)\big)
=\sum_{k=1}^p\Big\{g\big[\eta(y)+ku\big]-g\big[\eta(y)+(k-1)u\big]\Big\}.$$
Now it follows from (4.3) that
$$\big|f(\xi)-f(\eta)\big|\leq\sum_{k=1}^p\Big|g\big[\eta(y)+ku\big]-g\big[\eta(y)+(k-1)u\big]\Big|\leq pA|u|
=A\big|\xi(y)-\eta(y)\big|.\ \ \square$$
\end{Corollary}

For $z=(z_1,\cdots,z_m)\in\mathbb C^m$ let
$|z|=\sqrt{|z_1|^2+\cdots+|z_m|^2}$ and $\mathscr
K_{M,\varepsilon}(\mathbb C^m,\mathbb C)=\big\{g\in\mathbb
C^{\mathbb C^m}:g(0)=0$; for $z,u\in\mathbb C^m$ with
$|u|\leq\varepsilon$ and $t\in\mathbb C$ with $|t|\leq1$,
$g(z+tu)=g(z)+sg(u)$ where $|s|\leq M|t|\big\}$. For $g\in\mathscr
K_{M,\varepsilon}(\mathbb C^m,\mathbb C)$ and $1\leq j\leq m$ define
$g_j:\mathbb C\rightarrow\mathbb C$ by
$g_j(w)=g((0,\cdots,0,\stackrel{(j)}{w},0,\cdots,0)),\
\forall\,w\in\mathbb C$, then $g_j\in\mathscr
K_{M,\varepsilon}(\mathbb C,\mathbb C)$.

If $k\in\mathbb N$ then $\big\{\mbox{multi-index
}\alpha:|\alpha|\leq k\big\}$ is a finite set
$\big\{\alpha_1,\alpha_2,\cdots,\alpha_{m_k}\big\}$ which is
lexicographically ordered such that $\alpha_1=(0,\cdots,0)$,
$\alpha_2=(0,\cdots,0,1)$, $\cdots$, $\alpha_{m_k}=(k,0,\cdots,0)$,
and we can write
$(z_{\alpha_1},z_{\alpha_2},\cdots,z_{\alpha_{m_k}})=(z_\alpha)_{|\alpha|\leq
k}$ in $\mathbb C^{m_k}$.

\begin{Theorem}
$M\geq1$, $\gamma(t)=Mt$ for $t\in\mathbb C$. If $f\in
C^\infty(\Omega)^{[\gamma,V]}$ is of order $k$ and $supp\,f=\{y\}$,
then there exist $\varepsilon>0$ and $g\in\mathscr
K_{M,\varepsilon}(\mathbb C^{m_k},\mathbb C)$ such that
$$f(\xi)=g\Big(\big(\partial^\alpha\xi(y)\big)_{|\alpha|\leq k}\Big),\ \ \forall\,\xi\in C^\infty(\Omega).$$

Proof. Letting $\eta_\alpha=(\bullet-y)^\alpha/(\alpha!)$ for
$|\alpha|\leq k$, by Cor. 4.3 we have
$$f(\xi)=f\Big(\sum_{|\alpha|\leq k}\partial^\alpha\xi(y)\eta_\alpha\Big),\ \ \forall\,\xi\in C^\infty(\Omega).$$
Pick a $U\in\mathcal N(C^\infty(\Omega))$ for which
$\stackrel{(m_k)}{\overbrace{U+U+\cdots+U}}\subset V$. There is an
$\varepsilon>0$ such that $u\eta_\alpha\in U$ when $u\in\mathbb C$
with $|u|\leq\varepsilon$ and $|\alpha|\leq k$. Hence
$$\sum_{|\alpha|\leq k}u_\alpha\eta_\alpha\in\stackrel{m_k}{\overbrace{U+U\cdots+U}}\subset V,\ \
\forall\,(u_\alpha)_{|\alpha|\leq k}\in \mathbb C^{m_k}\mbox{ with
}|(u_\alpha)_{|\alpha|\leq k}|\leq\varepsilon.$$

Define $g:\mathbb C^{m_k}\rightarrow\mathbb C$ by
$$g\Big(\big(z_\alpha\big)_{|\alpha|\leq k}\Big)=f\Big(\sum_{|\alpha|\le
k}z_\alpha\eta_\alpha\Big),\ \ \forall\,(z_\alpha)_{|\alpha|\leq
k}\in\mathbb C^{m_k}.$$ Then for $(z_\alpha)_{|\alpha|\leq k}$,
$(u_\alpha)_{|\alpha|\leq k}\in\mathbb C^{m_k}$ with
$|(u_\alpha)_{|\alpha|\leq k}|\leq\varepsilon$ and $t\in\mathbb C$
with $|t|\leq1$, $\sum_{|\alpha|\leq k}u_\alpha\eta_\alpha\in V$ and
\begin{align*}
g\Big(\big(z_\alpha\big)_{|\alpha|\leq
k}+t\big(u_\alpha\big)_{|\alpha|\leq
k}\Big)&=g\Big(\big(z_\alpha+tu_\alpha\big)_{|\alpha|\leq
k}\Big)=f\Big(\sum_{|\alpha|\leq
k}\big(z_\alpha+tu_\alpha\big)\eta_\alpha\Big)\\
&=f\Big(\sum_{|\alpha|\leq k}z_\alpha\eta_\alpha+t\sum_{|\alpha|\leq
k}u_\alpha\eta_\alpha\Big)\\
&=f\Big(\sum_{|\alpha|\leq
k}z_\alpha\eta_\alpha\Big)+sf\Big(\sum_{|\alpha|\leq
k}u_\alpha\eta_\alpha\Big)\\
&=g\Big(\big(z_\alpha\big)_{|\alpha|\leq
k}\Big)+sg\Big(\big(u_\alpha\big)_{|\alpha|\leq k}\Big),
\end{align*}
where $|s|\leq|\gamma(t)|\leq M|t|$.

Thus, $g\in\mathscr K_{M,\varepsilon}(\mathbb C^{m_k},\mathbb C)$
and
$$f(\xi)=f\Big(\sum_{|\alpha|\leq k}\partial^\alpha\xi(y)\eta_\alpha\Big)=
g\Big(\big(\partial^\alpha\xi(y)\big)_{|\alpha|\leq k}\Big),\ \
\forall\,\xi\in C^\infty(\Omega).\ \ \square$$
\end{Theorem}

\begin{Theorem}
$M\geq1$, $\gamma(t)=Mt$ for $t\in\mathbb C$. If $f\in
C^\infty(\Omega)^{[\gamma,V]}$ is of order $k$ and $supp\,f=\{y\}$,
then there exists 
$\big\{g_\alpha:\alpha\mbox{ is a multi-index},\ |\alpha|\leq
k\big\}\subset\mathscr K_{M,\varepsilon}(\mathbb C,\mathbb C)$ such
that
$$f\big[z(\bullet-y)^\alpha/(\alpha!)\big]=g_\alpha(z),\ \ \forall\,|\alpha|\leq k,\ z\in\mathbb C,$$
$$f(\xi)=M\sum_{|\alpha|\leq
k}a_{\alpha,\xi}\frac{g_\alpha(\varepsilon)}{\varepsilon}\partial^\alpha\xi(y),\
\ \forall\,\xi\in C^\infty(\Omega),$$ where all
$|a_{\alpha,\xi}|\leq1$.

Proof. Letting $\eta_\alpha=(\bullet-y)^\alpha/\alpha!$ for
$|\alpha|\leq k$, by Cor. 4.3 we have
$$f(\xi)=f\Big(\sum_{|\alpha|\leq k}\partial^\alpha\xi(y)\eta_\alpha\Big),\ \ \forall\,\xi\in C^\infty(\Omega).$$

There is an $\varepsilon>0$ such that $u\eta_\alpha\in V$ for all
$u\in\big\{z\in\mathbb C:|z|\leq\varepsilon\big\}$ and $|\alpha|\leq
k$. We write $\big\{\alpha:|\alpha|\leq
k\big\}=\big\{\alpha_1,\alpha_2,\cdots,\alpha_m\big\}$.

Let $\xi\in C^\infty(\Omega)$ and pick a $p\in\mathbb N$ such that
$|\frac{\partial^\alpha\xi(y)}{p\varepsilon}|<1$ when $|\alpha|\leq
k$. Then
\begin{align*}
f(\xi)&=f\Big(\sum_{j=1}^m\partial^{\alpha_j}\xi(y)\eta_{\alpha_j}\Big)=
f\Big(\sum_{j=1}^{m-1}\partial^{\alpha_j}\xi(y)\eta_{\alpha_j}+
p\frac{\partial^{\alpha_m}\xi(y)}{p\varepsilon}\varepsilon\eta_{\alpha_m}\Big)\\
&=f\Big(\sum_{j=1}^{m-1}\partial^{\alpha_j}\xi(y)\eta_{\alpha_j}\Big)
+\big(\sum_{v=1}^ps_v\big)f(\varepsilon\eta_{\alpha_{m}}),
\end{align*}
where each $|s_v|\leq
M|\frac{\partial^{\alpha_m}\xi(y)}{p\varepsilon}|$ so
$|\sum_{v=1}^ps_v|\leq
pM|\frac{\partial^{\alpha_m}\xi(y)}{p\varepsilon}|
=\frac{M}{\varepsilon}|\partial^{\alpha_m}\xi(y)|$.

Define $g_{\alpha_m}:\mathbb C\rightarrow\mathbb C$ by
$g_{\alpha_m}(z)=f(z\eta_{\alpha_m}),\ \forall\,z\in\mathbb C$. For
$z,u,t\in\mathbb C$ with $|u|\leq\varepsilon$ and $|t|\leq1$,
$u\eta_{\alpha_m}\in V$ and
$g_{\alpha_m}(z+tu)=f(z\eta_{\alpha_m}+tu\eta_{\alpha_m})=
f(z\eta_{\alpha_m})+sf(u\eta_{\alpha_m})=g_{\alpha_m}(z)+sg_{\alpha_m}(u)$
where $|s|\leq|\gamma(t)|=M|t|$. Thus, $g_{\alpha_m}\in\mathscr
K_{M,\varepsilon}(\mathbb C,\mathbb C)$ and
$$\Big|\big(\sum_{v=1}^ps_v\big)f(\varepsilon\eta_{\alpha_m})\Big|
=\Big|\big(\sum_{v=1}^ps_v\big)g_{\alpha_m}(\varepsilon)\Big|\leq
M\Big|\frac{g_{\alpha_m}(\varepsilon)}{\varepsilon}\partial^{\alpha_m}\xi(y)\Big|.$$
Hence there is an $a_{\alpha_m,\xi}\in\mathbb C$ such that
$|a_{\alpha_m,\xi}|\leq1$ and
$(\sum_{v=1}^ps_v)f(\varepsilon\eta_{\alpha_m})=
a_{\alpha_m,\xi}M\frac{g_{\alpha_m}(\varepsilon)}{\varepsilon}\partial^{\alpha_m}\xi(y)$.
In this way, we have
\begin{align*}
f(\xi)&=f\Big(\sum_{j=1}^{m-1}\partial^{\alpha_j}\xi(y)\eta_{\alpha_j}\Big)+
a_{\alpha_m,\xi}M\frac{g_{\alpha_m}(\varepsilon)}{\varepsilon}\partial^{\alpha_m}\xi(y)\\
&=f\Big(\sum_{j=1}^{m-2}\partial^{\alpha_j}\xi(y)\eta_{\alpha_j}\Big)+
a_{\alpha_{m-1},\xi}M\frac{g_{\alpha_{m-1}}(\varepsilon)}{\varepsilon}\partial^{\alpha_{m-1}}\xi(y)
+a_{\alpha_m,\xi}M\frac{g_{\alpha_m}(\varepsilon)}{\varepsilon}\partial^{\alpha_m}\xi(y)\\
&\qquad\qquad\cdots\\
&=\sum_{j=1}^ma_{\alpha_j,\xi}M\frac{g_{\alpha_j}(\varepsilon)}{\varepsilon}\partial^{\alpha_j}\xi(y)\\
&=M\sum_{|\alpha|\leq
k}a_{\alpha,\xi}\frac{g_{\alpha}(\varepsilon)}{\varepsilon}\partial^{\alpha}\xi(y),
\end{align*}
where all $|a_{\alpha,\xi}|\leq1$.\ \ $\square$
\end{Theorem}

It is similar to Cor. 4.5 that we have

\begin{Corollary}
$M\geq1$, $\gamma(t)=Mt$ for $t\in\mathbb C$. If $f\in
C^\infty(\Omega)^{[\gamma,V]}$ is of order $k$ and $supp\,f=\{y\}$,
then $f$ has the following properties.

(1) If $|\alpha|\leq k$ and
$f[z_0(\bullet-y)^\alpha/(\alpha!)]\neq0$ for some $z_0\in\mathbb
C$, then there exists an $\varepsilon>0$ such that
$f[z_0(\bullet-y)^\alpha/(\alpha!)]\neq0$ when
$0<|z|\leq\varepsilon$, that is, the equation
$$f\big[z_0(\bullet-y)^\alpha/(\alpha!)\big]=0,\ \ |z|\leq\varepsilon$$
has the unique solution $z=0$. Hence if $|\alpha|\leq k$ and there
is $\{z_v\}\subset\mathbb C$ such that each $z_v\neq0$,
$z_v\rightarrow0$ and each $f[z_v(\bullet-y)^\alpha/(\alpha!)]=0$,
then $f[z(\bullet-y)^\alpha/(\alpha!)]=0$ for all $z\in\mathbb C$.

(2) If $|\alpha|\leq k$, then $f[z(\bullet-y)^\alpha/(\alpha!)]$ is
Lipschitz, that is, there is an $A_\alpha>0$ such that
$$\Big|f\big[z(\bullet-y)^\alpha/(\alpha!)\big]-f\big[u(\bullet-y)^\alpha/(\alpha!)\big]\Big|
\leq A_\alpha|z-u|,\ \ \forall\,z,u\in\mathbb C.$$ In particular, we
have
$$\Big|f\big[\partial^\alpha\xi(y)(\bullet-y)^\alpha/(\alpha!)\big]
-f\big[\partial^\alpha\eta(y)(\bullet-y)^\alpha/(\alpha!)\big]\Big|\leq
A_\alpha\big|\partial^\alpha\xi(y)-\partial^\alpha\eta(y)\big|,\
\forall\,\xi,\eta\in C^\infty(\Omega).$$
\end{Corollary}

\end{document}